\documentclass[12pt]{article}
\usepackage{amsmath}
\usepackage{amssymb}
\usepackage{algorithm}
\usepackage{algorithmic}
\usepackage[dvips]{graphicx}
\usepackage{subfigure}

\def\eps{\varepsilon}
\def\bbox{\quad\hbox{\vrule \vbox{\hrule \vskip2pt \hbox{\hskip2pt
\vbox{\hsize=1pt}\hskip2pt} \vskip2pt\hrule}\vrule}}
\def\lessim{\ \lower4pt\hbox{$
\buildrel{\displaystyle <}\over\sim$}\ }
\def\gessim{\ \lower4pt\hbox{$\buildrel{\displaystyle >}
\over\sim$}\ }

\def\HH{{\cal H}}
\def\eps{{\varepsilon}}

\newcommand{\Reals}{\mathbb{R}}

\newtheorem{lemma}{Lemma}
\newtheorem{theorem}{Theorem}
\newtheorem{corollary}{Corollary}
\makeatletter
\@addtoreset{equation}{section}

\makeatother

\font\tencmmib=cmmib10 \skewchar\tencmmib '60
\newfam\cmmibfam
\textfont\cmmibfam=\tencmmib

\def\bbox{\quad\hbox{\vrule \vbox{\hrule \vskip2pt \hbox{\hskip2pt
\vbox{\hsize=1pt}\hskip2pt} \vskip2pt\hrule}\vrule}}
\def\lessim{\ \lower4pt\hbox{$
\buildrel{\displaystyle <}\over\sim$}\ }
\def\gessim{\ \lower4pt\hbox{$\buildrel{\displaystyle >}
\over\sim$}\ }

\def\eps{\varepsilon}

\def\go0{\to 0}

\def\HH{{\cal H} }

\def\leftitem#1{\item{\hbox to\parindent{\enspace#1\hfill}}}

\def\mdsk{\medskip}

\def\qed{\hfill\break\rightline{$\bbox$}}

\def\sg{\sigma}

\def\sg2{\sigma^2}

\def\__{_{\infty}}
\begin{document}

\title{Bounding the generalization error of convex combinations
of classifiers: balancing the dimensionality and the margins}

\author{
Vladimir Koltchinskii
\thanks{Partially supported by NSA Grant MDA904-99-1-0031},
\ \ Dmitriy Panchenko
\thanks{Partially supported by UNM Office of Graduate Studies RPT Grant}
\\ Department of Mathematics and Statistics
\\ The University of New Mexico
\\ and Fernando Lozano
\\ Department of Electrical and Computer Engineering
\\ The University of New Mexico
}

\maketitle
\begin{abstract}
A problem of bounding the generalization error of a classifier
$f\in {\rm conv}({\cal H}),$ where ${\cal H}$ is a "base" class
of functions (classifiers), is considered. This problem frequently
occurs in computer learning, where efficient algorithms of combining
simple classifiers into a complex one (such as boosting and bagging)
have attracted a lot of attention. Using Talagrand's concentration
inequalities for empirical processes, we obtain new sharper bounds
on the generalization error of combined classifiers that take into
account both the empirical distribution of ``classification margins''
and an "approximate dimension" of the classifiers and study
the performance of these bounds in several experiments with
learning algorithms.

\end{abstract}

\vskip 5mm

\hfill\break
{\it 1991 AMS subject classification}: primary 62G05,
secondary 62G20, 60F15
\hfill\break
{\it Keywords and phrases}: generalization error,
combined classifier, margin, approximate dimension,
empirical process, Rademacher process,
random entropies, concentration inequalities,
boosting, bagging
\hfill\break
{\it Abbreviated Title}: Dimensionality and Margins

\vfill\break

\section{Introduction}

Let $(X_1,Y_1),\dots ,(X_n,Y_n)$ be a sample of $n$ labeled
training examples that are independent identically distributed
copies of a random couple $(X,Y),$ $X$ being an ``instance'' in
a measurable space $S$ and
$Y$ being a ``label'' taking values in $\{-1,1\}.$ Let $P$ denote
the distribution of the couple $(X,Y).$
Given a measurable function $f$ from $S$ into ${\mathbb R},$ we use
${\rm sign}(f(x))$ as a predictor of the unknown label
of an instance $x\in S.$ We will call $f$ a classifier
of the examples from $S.$
The quantity ${\mathbb P}\{Yf(X)\leq 0\}={\mathbb P}\{(x,y):yf(x)\leq 0\}$
is called \it the generalization error \rm of the classifier
$f.$ The goal of learning (classification) is, given a set of
training examples, to find a classifier $f$ with a small generalization
error.

Some of the important recent advances in statistical learning
theory are related to the development of complex classifiers
that are combinations of simpler ones. In so called \it voting
methods \rm of combining classifiers (such as boosting, bagging,
etc.) a complex classifier produced by a learning algorithm
is a convex combination of simpler classifiers from the base
class.

Let ${\cal H}$ be a class of functions
from $S$ into ${\mathbb R}$ (base classifiers) and let
${\cal F}:={\rm conv}({\cal H})$ denote the symmetric
convex hull of ${\cal H}:$
$$
{\rm conv}({\cal H}):=
\Bigl\{\sum_{i=1}^N \lambda_i h_i: N\geq 1, \lambda_i\in {\mathbb R},
\sum_{i=1}^N |\lambda_i|\leq 1,\ h_i\in {\cal H}\Bigr\}.
$$

Our main goal in this paper is to develop new probabilistic
upper bounds on the generalization error of a classifier
$f$ from the symmetric convex hull ${\cal F}={\rm conv}({\cal H})$
of the base class. The well known approach to such a problem,
developed in pathbreaking works of Vapnik and Chervonenkis
(see \cite{VAP98} and references therein), is based on
an easy bound
$$
P\{(x,y):yf(x)\leq 0\}\leq P_n\{(x,y):yf(x)\leq 0\}+
\sup_{C\in {\cal C}}[P(C)-P_n(C)],
$$
where $P_n$ is the empirical distribution of the training
examples, i.e. for any set $C\subset S\times \{-1,1\},$
$P_n(C)$ is the frequency of training examples in the set $C,$
$$
{\cal C}:= \Bigl\{\{(x,y):yf(x)\leq 0\}:f\in {\cal F}\Bigr\},
$$
and on further bounding the uniform (over the class ${\cal C}$)
deviation of the empirical distribution $P_n$ from the true
distribution $P.$ The methods that are used to solve this
problem belong to the
theory of empirical processes and the crucial role is played
by the VC-dimension of the class ${\cal C},$ or by more
sophisticated entropy characteristics of the class.
For instance, if $m^{{\cal C}}(n)$ denotes the maximal
number of subsets obtainable by intersecting a sample
of size $n$ with the class ${\cal C}$ (the so called shattering
number),
%and $\hat f$ is a classifier that minimizes the training
%error $P_n\{(x,y):yf(x)\leq 0\}$ on the class ${\cal F},$
then
the following bound holds
(see \cite{DVGL}, Theorem 12.6) for all $\eps >0$
$$
{\mathbb P}\Bigl\{
P\{(x,y):yf(x)\leq 0\}\geq P_n\{(x,y):yf(x)\leq 0\}+\eps\Bigr\}\leq
8m^{{\cal C}}(n)e^{-n\eps^2/32}.
$$
It follows from this bound that the training error measures the generalization
error of a classifier $f\in {\cal F}$ with the accuracy
$O\Bigl(\sqrt{\frac{V({\cal C})\log n}{n}}\Bigr),$ where $V({\cal C})$
is the VC-dimension of the class ${\cal C}.$
In the so called zero-error case, when there exists
a classifier $\hat f\in {\cal F}$ with zero training error,
we even have the bound
(see \cite{DVGL}, Theorem 12.7):
$$
{\mathbb P}\Bigl\{P\{(x,y):y\hat f(x)\leq 0\}\geq \eps\Bigr\}\leq
2m^{{\cal C}}(2n)2^{-n\eps /2},
$$
which implies that the generalization error of the classifier
$\hat f$ is of the order $O\Bigl(\frac{V({\cal C})\log n}{n}\Bigr).$
The above bounds, however, do not apply directly to the case
of the class ${\cal F}={\rm conv}({\cal H}),$ which is of interest in
applications to bounding the generalization error of the
voting methods, since in this case typically $V({\cal C})=+\infty .$
Even when one deals with a finite number of base classifiers in
a convex combination (which is the case, say, with boosting
after finite number of rounds), the VC-dimensions of the classes
involved are becoming rather large, so the above bounds do not
explain the generalization ability of boosting and other voting
methods observed in numerous experiments. This motivated
Bartlett \cite{Bartlett}, Schapire, Freund, Bartlett and Lee \cite{SFBL}
(see also \cite{Anthony:Bartlett})
to develop a new class of upper bounds
on generalization error of a convex combination of classifiers,
expressed in terms of empirical distribution of margins
(the role of classification margins in improving the
generalization ability
of learning machines was clear in earlier work on support vector
machines as well, see \cite{Cortes:Vapnik}.
The margin
of a classifier $f$ on a training example $(X,Y)$ is defined
as the product $Yf(X).$ Schapire, Freund, Bartlett and Lee \cite{SFBL}
showed that for a given $\alpha \in (0,1)$ with probability
at least $1-\alpha$ for all $f\in {\rm conv}({\cal H})$
$$
P\{(x,y):yf(x)\leq 0\}\leq
\inf_{\delta}\Bigl[P_n\{(x,y):yf(x)\leq \delta\}+
\frac{C}{\sqrt{n}}\Bigl(\frac{V({\cal H})
\log^2(\frac{n}{V({\cal H})})}{\delta^2}
+\log(1/\alpha)
\Bigr)^{1/2}
\Bigr].
$$
Choosing in the above bound the value of $\delta = \hat \delta (f)$ that solves
the equation
$$
\delta P_n\{(x,y):yf(x)\leq \delta\} = \sqrt{\frac{V({\cal H})}{n}}
$$
(which is nearly an optimal choice), one gets (ignoring the logarithmic
factors) the generalization error of a classifier $f$
from the convex hull of the order
$$O\Bigl(\frac{1}{\hat \delta(f)}\sqrt{\frac{V({\cal H})}{n}}\Bigr).$$
Koltchinskii and Panchenko \cite{KP}, using the methods of the theory
of Empirical, Gaussian and Rademacher Processes (concentration
inequalities, symmetrization, comparison inequalities)
generalized and refined this type of bounds. They also suggested
a way to improve these bounds under certain assumptions on the growth
of random entropies of a class ${\cal F}$ to which the classifier
belongs. The new bounds are based on the notion of $\gamma$-margin
of the classifier, introduced in their paper. The $\gamma$-margins
are defined for $\gamma\in (0,1)$ (see the definitions in Section 2
below), the value of $\gamma=1$ roughly corresponds to the case
studied in \cite{SFBL}. The quality
of the bound improves as $\gamma$ decreases to $0.$ However, the
bounds of this type are proved to hold for the values of
$\gamma \geq 2\alpha/(2+\alpha),$ where $\alpha\in (0,2)$ is
the growth exponent of the random entropy of the class ${\cal F}.$
In the case of ${\cal F}:={\rm conv}({\cal H}),$ where ${\cal H}$
is a VC-class with VC-dimension $V({\cal H}),$ this leads to
the values of $\alpha ={2(V({\cal H})-1)/V({\cal H})}<2,$ which
allows one to use $\gamma$-margins with $\gamma<1$ (but it is
going to be rather close to $1$ unless the VC-dimension is very
small). The experiments of Koltchinskii, Panchenko and Lozano \cite{KPL}
showed that, in the case of the classifiers obtained in consecutive rounds
of boosting, the bounds on the generalization error
in terms of $\gamma$-margins hold even for much smaller values of
$\gamma .$ This allows one to conjecture that such classifiers
belong, in fact, to a class ${\cal F}\subset {\rm conv}({\cal H})$
whose entropy might be much smaller than the entropy of the whole
convex hull. The problem, though, is that it is practically
impossible to identify such a class prior to experiments, leaving
the question of how to choose the values of $\gamma$ for which the
bounds hold open. In this paper, we develop a new approach
to this problem. Namely, we suggest an adaptive bound on the generalization
error of a convex combination of classifiers from a base class
that is based on the one hand on the margins of the combined classifiers
and on the other hand on their \it approximate dimensions \rm
(the numbers of ``large enough'' coefficients in the convex combinations).
This adaptive bound ``captures'' the size of the entropy
of a subset of the convex hull to which the classifier actually
belongs.

The results are formulated precisely in Section 2. The proofs
that heavily rely upon Talagrand's concentration and deviation
inequalities for empirical processes are given in section 3.
Section 4 includes the results of several experiments with
existing learning algorithms (such as boosting and bagging)
for which we computed the bounds on the learning curves that
follow from our results. We also discuss here some approaches
to combining classifiers that attempt to minimize the margin
cost function keeping the dimension of the classifier small.

\section{Empirical margins and approximate dimensions: main results}

Let $(S,{\cal A})$ be a measurable space and let ${\cal F}$ be a
class of measurable functions on $(S,{\cal A})$. In this section,
in order to shorten the notations, we suppress the labeles. If
one wants to apply the results in the setting of the Introduction,
one has to consider instead of $S$ the space $S\times \{-1,1\}$
and instead of a function $f$ on $S,$ a function $(x,y)\mapsto yf(x)$
on $S\times \{-1,1\}.$ The results can be also used in the case of
multiclass problems (see Section 5 in \cite{KP}).
In what follows $P$ denotes a probability measure on $(S,{\cal A}),$
$\{X_n\}$ is a sequence of i.i.d. random variables, defined on a
probability space $(\Omega,\Sigma,{\mathbb P})$ and taking values
in $(S,{\cal A})$ with distribution $P,$ $P_n$ denote the empirical
measure based on the sample $(X_1,\dots,X_n):$
$$
P_n(A):=n^{-1}\sum_{i=1}^n I_A(X_i),\ A\subset S.
$$

We start with extending the bounds on generalization error,
obtained by Koltchinskii and Panchenko \cite{KP} in terms of so called
$\gamma$-margins.

Below we give a definition of what we call
$\psi -$bounds that will play a major role
in bounding the generalization error of classifiers.
These quantities depend on a function $\psi$ that will
characterize the complexity of the class ${\cal F},$
and therefore determine the quality of the bounds.

Let $\psi $ be a concave nondecreasing function
on $[0,+\infty)$ with $\psi (0)=0.$
%For a fixed $\delta>0,$
%denote by $\eps_n^{\psi}(\delta)$ the smallest solution of the equation
For a fixed $\eps>0,$
denote by $\delta_n^{\psi}(\eps)$ the largest solution of the equation
\begin{equation}
\eps=\frac{1}{\delta\sqrt{n}}\psi(\delta\sqrt{\eps})
\label{Equa}
\end{equation}
(if $\psi$ is strictly concave,
the solution of the equation (\ref{Equa})
is unique).
Clearly, for a concave $\psi$ the function
$\varphi (x)\equiv \frac{\psi(x)}{x}$ is nonincreasing.
Therefore, it is easy to see that
$$
\delta_n^{\psi}(\eps)=\frac{\varphi^{-1}(\sqrt{\eps n})}{\sqrt{\eps}}.
$$
Given a function $f$ and $t>0,$ define the following quantity
$$
\eps_n^{\psi}(f;t):=\inf\Bigl\{\eps \geq \frac{t\bigvee 2\log n}{n}:
P\{f\leq \delta_n^{\psi}(\eps)\}\leq \eps\Bigr\}
$$
and its empirical version
$$
\hat \eps_n^{\psi}(f;t):=
\inf\Bigl\{\eps \geq \frac{t\bigvee 2\log n}{n}:
P_n\{f\leq \delta_n^{\psi}(\eps)\}\leq \eps\Bigr\}
$$
Since for all $\eps>0,$ $\delta_n^{\psi}(\eps)\geq 0,$
it immediately follows from the definition that for all $f\in {\cal F}$
$$
P\{f\leq 0\}\leq
\inf\{
P\{f\leq \delta_n^{\psi}(\eps)\}:
\eps\geq \eps_n^{\psi}(f;t)
\}\leq
\eps_n^{\psi}(f;t).
$$
We will call $\eps_n^{\psi}(f;t)$ and $\hat \eps_n^{\psi}(f;t)$
\it the $\psi$-bound \rm and \it the empirical $\psi$-bound \rm
of the classifier $f,$ respectively.
We show below that under a proper assumption on the
random entropy of the class ${\cal F},$ with a high
probability the empirical $\psi$-bounds $\hat \eps_n^{\psi}(f;t)$
are, for all the functions from the class, within
a multiplicative constant from the true $\psi$-bounds $\eps_n^{\psi}(f;t).$
This allows one to replace $\eps_n^{\psi}(f;t)$
 in the above bound on $P\{f\leq 0\}$
by $\hat \eps_n^{\psi}(f;t)$ (which gives in applications
a bound on the generalization errors of classifiers).

Given a metric space $(T,d),$ we denote $H_d(T;\eps)$
the $\eps$-entropy of $T$ with respect to $d,$ i.e.
$$
H_d(T;\eps):=\log N_d(T;\eps),
$$
where $N_d(T;\eps)$ is the minimal number of balls of
radius $\eps$ covering $T.$
If $Q$ is a probability measure on $(S;{\cal A}),$
$d_{Q,2}$ will denote the metric of the space $L_2(S;dQ):$
$
d_{Q,2}(f;g):=(Q|f-g|^2)^{1/2}.
$

\begin{theorem}\label{theorem:psi}
 Let $\psi$ be a concave nondecreasing function
on $[0,+\infty)$ with $\psi(0)=0.$
%$\psi(x)\geq 3x\sqrt{\log\frac{e}{x}},\ x\in (0,1).$
Suppose the following bound on Dudley's entropy integral holds
with some $D_n>0:$
\begin{equation}
\int\limits_{0}^{x}H_{d_{P_n, 2}}^{1/2}({\cal F},u)du\leq D_n\psi(x),\
x>0 \ {\rm a.s.} \label{entcond}
\end{equation}
where $D_n=D_n(X_1,\ldots,X_n)$ is a function of training
examples such that ${\mathbb E}D_n< \infty.$
Then there exist absolute constants $A, B>0$
such that for $\bar A:=A(1+{\mathbb E}D_n)^2$ and for all $t>0$
\begin{eqnarray}
&&
{\mathbb P}\Bigl\{\forall f\in {\cal F}:
{\bar A}^{-1}\hat \eps_n^{\psi}(f;t)\leq
\eps_n^{\psi}(f;t)\leq
\bar A \hat \eps_n^{\psi}(f;t)
\Bigr\}
\nonumber
\\
&&
\geq 1-B\log_2\log_2 {\frac{n}{t\bigvee 2\log n}}
\exp\Bigl\{-\bigl(\frac{t}{2}\bigvee \log n\bigr)\Bigr\}.
\label{psibound}
\end{eqnarray}
\end{theorem}

The following corollary is immediate.

\begin{corollary}
\label{corollary:psi}
Under the conditions of Theorem 1
there exist numerical constants $A,B>0$
such that for $\bar A:=A(1+{\mathbb E}D_n)^2$
and for all $t>0$
\begin{equation}
{\mathbb P}\Bigl\{\exists f\in {\cal F}:
P\{f\leq 0\}\geq \bar A\hat \eps_n^{\psi}(f;t)
\Bigr\}
\leq
B\log_2\log_2 {\frac{n}{t\bigvee 2\log n}}
\exp\Bigl\{-\bigl(\frac{t}{2}\bigvee \log n\bigr)\Bigr\}.
\label{ge1}
\end{equation}
\end{corollary}

\mdsk

{\bf Example 1}. Let $\alpha \in (0,2)$ and
$\psi (x)\equiv x^{1-\alpha/2}.$
Let $\gamma :=\frac{2\alpha}{\alpha+2}.$
Koltchinskii and Panchenko \cite{KP}
defined $\gamma$-margins of a function $f$ as follows:
$$
\delta_n(\gamma ;f):=\sup\Bigl\{\delta \in (0,1):
\delta^{\gamma}P\{f\leq \delta\}\leq n^{-1+\frac{\gamma}{2}}\Bigr\},
$$
$$
\hat \delta_n(\gamma ;f):=\sup\Bigl\{\delta \in (0,1):
\delta^{\gamma}P_n\{f\leq \delta\}\leq n^{-1+\frac{\gamma}{2}}\Bigr\}.
$$
An easy computation shows that
$$
\eps_n^{\psi}(f;n^{\gamma/2})=\frac{1}{n^{1-\gamma/2}\delta_n(\gamma;f)^{\gamma}}.
$$
Corollary 1 immediately implies that if for some
$\alpha \in (0,2)$ and $D_n>0,$
${\mathbb E}D_n < \infty$
$$
H_{d_{n,2}}({\cal F};u)\leq D_n^2 u^{-\alpha},\ u>0\ {\rm a.s.},
$$
then for any $\gamma \geq \frac{2\alpha}{\alpha+2}$ there
exist constants $A, B>0$ such that for $\bar A:=A(1+{\mathbb E}D_n)^2$
\begin{equation}
\label{gamma1}
{\mathbb P}\Bigl\{\exists f\in {\cal F}:
P\{f\leq 0\}\geq \frac{\bar A}{n^{1-\gamma/2}\hat \delta_n(\gamma ;f)^{\gamma}}
\Bigr\}
\leq B\log_2\log_2 n\exp\Bigl\{-n^{\gamma/2}/2\Bigr\}
\end{equation}
(see also \cite{KP}).
It is easy to see that the quantity
\begin{equation}
\label{eq:gammabound}
\frac{1}{n^{1-\gamma/2}\hat \delta_n(\gamma ;f)^{\gamma}}
\end{equation}
in the above upper bound on the generalization
error \it becomes smaller \rm as $\gamma$ decreases from $1$ to $0.$
The Schapire-Freund-Bartlett-Lee type of bounds correspond
to the worst choice of $\gamma$ ($\gamma=1$). In the case when
${\cal F}$ is the symmetric convex hull of a VC-class ${\cal H}$
with VC-dimension $V({\cal H})$
the value of ${\alpha}$ is equal to
$\frac{2(V({\cal H})-1))}{V({\cal H})}<2$
that allows us to have $\gamma <1,$ improving the previously known
bound. In fact, Koltchinskii, Panchenko and Lozano \cite{KPL} computed
the empirical $\gamma$-margins of classifiers obtained in consecutive rounds
of boosting and observed that the bounds on their generalization error
in terms of $\gamma$-margins hold even for much smaller values of
$\gamma .$ This allows one to conjecture that such classifiers
belong, in fact, to a class ${\cal F}\subset {\rm conv}({\cal H})$
whose entropy might be much smaller than the entropy of the whole
convex hull.

{\bf Example 2}. Consider now the case of
$\psi (x)\equiv x\sqrt{\log\frac{e}{x}}$ for
$x\leq 1$ and $\psi (x)\equiv x$ for $x>1.$
Then, by a simple computation,
$$
\delta_n^{\psi}(\eps)=\frac{e^{1-n\eps}}{\sqrt{\eps}},\ \eps\geq n^{-1}.
$$
If we define
\begin{equation}
\label{eq:epsVC}
\hat \eps_n^{VC}(f;t):=
\inf\Bigl\{\eps \geq \frac{t\bigvee 2\log n}{n}:
P_n\{f\leq \frac{e^{1-n\eps}}{\sqrt{\eps}}\}\leq \eps\Bigr\},
\end{equation}
then under the condition
$$
H_{d_{P_n,2}}({\cal F};u)\leq D_n^2 \log\frac{1}{u}\bigvee 1,\ u>0\ {\rm a.s.},
$$
with some $D_n=D_n(X_1,\dots, X_n),\ {\mathbb E}D_n<+\infty $
(which holds, for instance, if ${\cal F}$ is a VC-subgraph class),
we get from Corollary 1 that with some numerical constants $A,B>0$
for all $t>0$
$$
{\mathbb P}\Bigl\{\exists f\in {\cal F}:
P\{f\leq 0\}\geq
\bar A\hat \eps_n^{VC}(f;t)
\Bigr\}
\leq
B\log_2\log_2 {\frac{n}{t\bigvee 2\log n}}
\exp\Bigl\{-\bigl(\frac{t}{2}\bigvee \log n\bigr)\Bigr\},
$$
where $\bar A:=A(1+{\mathbb E}D_n)^2.$

The proofs of Theorem 1 and Theorem 3 below are based on the
following generalization of one of the results of Koltchinskii
and Panchenko \cite{KP}
(that itself relies heavily on the concentration inequality
for empirical processes due to Talagrand).

Given a nondecreasing concave function $\psi$ on $[0,+\infty)$
with $\psi (0)=0$ and a fixed number $\delta >0,$ we denote
by $\eps_n^{\psi}(\delta)>0$ the smallest solution of the equation
(\ref{Equa}) with respect to $\eps.$

\begin{theorem} Suppose that condition (\ref{entcond}) holds
with some concave nondecreasing $\psi$ such that $\psi(0)=0.$
Then, for all $\delta >0$ and for all
$\eps\geq \eps_n^{\psi}(\delta)\vee \frac{2\log n}{n}$
the following bounds hold
$$
{\mathbb P}\Bigl\{\exists f\in {\cal F}\
P_n\{f\leq \delta\}\leq
\varepsilon\
{\rm and}\
P\{f\leq  \frac{\delta}{2}\}\geq \bar A
\varepsilon
\Bigr\}\leq
$$
$$
\leq B \log_2\log_2\varepsilon^{-1}
\exp\{-\frac{n\varepsilon}{2}\}.
$$
and
$$
{\mathbb P}\Bigl\{\exists f\in {\cal F}\
P\{f\leq \delta\}\leq
\varepsilon\
{\rm and}\
P_n\{f\leq  \frac{\delta}{2}\}\geq \bar A\varepsilon
\Bigr\}\leq
$$
$$
\leq B \log_2\log_2\varepsilon^{-1}
\exp\{-\frac{n\varepsilon}{2}\},
$$
where $\bar A=A(1+{\mathbb E}D_n)^2$
and $A, B$ are numerical constants.
\end{theorem}

There are two major problems with the margin type bounds,
given above.
First of all, the values of the constants involved in the bounds
are far from being optimal and are too large at the moment.
Their improvement is related to a hard problem of optimizing the
constants in Talagrand's concentration inequalities for empirical
and Rademacher processes, used in the proofs below.
However, in the case when ${\cal F}={\rm conv}({\cal H})$
the constants in question depend only on the base class ${\cal H}$
and this allows one to use the bounds to study the behavior of the
generalization error when the the number of rounds of learning
algorithms (such as boosting) increases. Another problem is
related to the fact that there is no much prior knowledge about
the subset of ${\rm conv}({\cal H})$ to which a classifier
created by boosting or another method of combining the
classifiers is going to belong. This makes one to use the
value of
\begin{equation}
\label{equation:max_gamma}
\gamma =\frac{2\alpha}{\alpha+2}=\frac{2(V({\cal H})-1)}{2V({\cal H})-1}
\end{equation}
which is very
close to $1$ unless the VC-dimension of the base is \it very \rm
small. Our major goal in the current paper is to address this
problem. We do this by proving a new upper bound on the
generalization error of a classifier that belongs to a convex
hull of a base class. The bound includes the sum of
two main terms. The first one is an
``approximate'' dimension" of the classifier (the number of
``large enough'' coefficients in the convex combination) divided
by the sample size.
The second term is related to the margins of the classifier.
Balancing these two terms allows us to get rather tight upper
bound that ``captures'' the size of the entropy of a class to which
the classifier actually belongs. It combines previously known
bounds in terms of VC-dimension (in zero-error case) and in terms
of margins and becomes close to one of these two bounds
in the extreme cases.

Let ${\cal H}$ be
a class of measurable functions from $(S,{\cal A})$ into
${\mathbb R}.$ Let ${\cal F}\subset {\rm conv}({\cal H}).$
For a function $f\in {\cal F}$
and a number $\Delta \in [0,1],$ we define
\it the approximate $\Delta$-dimension \rm
of $f$ as the integer number $d\geq 0$ such that there
exist $N\geq 1,$ functions $h_j\in {\cal H},\ j=1,\dots ,N$ and numbers
$\lambda_j\in {\mathbb R},\ j=1,\dots ,N$ satisfying the conditions
$f=\sum_{j=1}^N\lambda_j h_j,$ $\sum_{j=1}^N |\lambda_j|\leq 1$ and
$\sum_{j=d+1}^N |\lambda_j|\leq \Delta .$ The $\Delta$-dimension
of $f$ will be denoted by $d(f;\Delta).$
Note that this definition depends on the representation
$f=\sum\lambda_j h_j,$ and one is free to use any
but the choice that produces smaller $d(f;\Delta)$
is advantageous.

In what follows we assume that for some $V>0$ and $K>0$
and for all probability measures $Q$ on $(S;{\cal A})$
\begin{equation}
N_{d_{Q,2}}({\cal H};(Q H^2)^{\frac{1}{2}}\eps)\leq K\eps^{-V},\ \eps>0,
\label{envc}
\end{equation}
where $H$ is a measurable envelope of ${\cal H}.$
In particular, this condition holds if ${\cal H}$ is a VC-subgraph
class. This condition implies the bound on the entropy
$$
H_{d_{Q,2}}({\rm conv}({\cal H});
(Q H^2)^{\frac{1}{2}}\eps)\leq C\eps^{-2V/(V+2)},\, \eps>0,
$$
where $C:=C(K;V)$ (see \cite{VandW}).
One can easily compute in this case that
$$
\int\limits_{0}^{x}H_{d_{P_n, 2}}^{1/2}({\cal F},u)du\leq
\frac{1}{2}(V+2)C^{1/2}(P_n H^2)^{\frac{V}{2(V+2)}}x^{\frac{2}{V+2}}
,\, x>0 \, {\rm a.s.}
$$
and, therefore, condition (\ref{entcond}) of Theorem 1
is satisfied with $\psi(x)= x^{\frac{2}{V+2}}$
under the assumption $P H^2 <\infty.$
Below we will assume that
one of the two conditions holds:

\begin{enumerate}
\item Class $\cal H$ is uniformly bounded and
${\cal F}\subset {\rm conv}({\cal H})$

\item The envelope $H$ of the class $\cal H$
is $P-$square integrable and
$$
{\cal F}\subset \Bigl\{\sum_{i=1}^N \lambda_i h_i:\ N\geq 1,
h_i\in {\cal H}, \lambda_i\in {\mathbb R}, \sum_{j=1}^N |\lambda_j|=1\Bigr\}.
$$
\end{enumerate}

Note, that under the second condition
$\cal F$ consists only of proper
symmetric convex combinations.

Let $\alpha :=\frac{2V}{V+2}$
and
$\Delta_f=\{\Delta \in [0,1] : d(f;\Delta)\leq n\}.$
Define
\begin{equation}
\eps_n(f;\delta):=\inf_{\Delta \in \Delta_f}\Bigl[
\frac{d(f;\Delta)}{n}\Bigl(\log\frac{1}{\delta}+
\log{\frac{ne^2}{d(f;\Delta)}}\Bigr)+
\Bigl(\frac{\Delta}{\delta}\Bigr)^{\frac{2\alpha}{\alpha+2}}
n^{-\frac{2}{\alpha+2}}\Bigr]\bigvee \frac{2\log n}{n}.
\label{efd}
\end{equation}
Let
$$
\hat \delta_n(f):=\sup\Bigl\{\delta\in (0,1/2): P_n\{f\leq \delta\}
\leq \eps_n(f;\delta)\Bigr\}.
$$

\begin{theorem}
\label{theorem:delta}
Assume that one of the above conditions on the class
$\cal F$ holds.
Then there exist constants $A,B>0$ such that for all
$0<t<n^{\frac{\alpha}{2+\alpha}}$
the following bound holds
$$
{\mathbb P}\Bigl\{\exists f\in {\cal F}\
P\{f\leq  \frac{\hat \delta_n(f)}{4}\}\geq \
A \Bigl(\eps_n(f;\frac{\hat \delta_n(f)}{2})+\frac{t}{n}\Bigr)
\Bigr\}\leq B e^{-t/4}.
$$
\end{theorem}

{\bf Example 3.}
If ${\cal F}\subset {\rm conv}({\cal H})$ is a class of functions such that
for some $\beta>0$
\begin{equation}
\sup_{f\in {\cal F}} d(f;\Delta)=O(\Delta^{-\beta}),
\label{plnmw}
\end{equation}
then with ``high probability'' for any classifier $f\in {\cal F}$
the upper bound on its generalization error becomes
of the order
$$
\frac{1}{n^{1-\gamma \beta/2(\gamma+\beta)}
\hat \delta_n(f)^{\gamma \beta/(\gamma +\beta)}},
$$
(which, of course, improves a more general bound in terms of
$\gamma$-margins; the general bound corresponds
to the case $\beta =+\infty$).
The condition (\ref{plnmw}) means that the weights of the
convex combination decrease polynomially fast, namely,
$|\lambda_j|= O(j^{-\alpha}),$ $\alpha=1+\beta^{-1}.$
The case of exponential decrease of the
weights is described by the condition
\begin{equation}
\sup_{f\in {\cal F}} d(f;\Delta)=O(\log\frac{1}{\Delta}).
\label{expw}
\end{equation}
In this case the upper bound becomes of the order
$
\frac{1}{n}\log^2 \frac{n}{\hat \delta_n(f)}.
$

\mdsk

\section{Proofs of the main results}

{\bf Proof of Theorem 1}. We use the first bound of Theorem 2.
The condition $\eps \geq \eps_n^{\psi}(\delta)$ is equivalent
to the condition $\delta \geq \delta_n^{\psi}(\eps).$ Thus,
we can use this bound for $\delta=\delta_n^{\psi}(\eps)$
and $\eps\geq (2\log n)/n.$ We get
$$
{\mathbb P}\Bigl\{\exists f\in {\cal F}\
P_n\{f\leq \delta_n^{\psi}(\eps)\}\leq
\varepsilon\
{\rm and}\
P\{f\leq  \frac{\delta_n^{\psi}(\eps)}{2}\}\geq \bar A
\varepsilon
\Bigr\}\leq
B \log_2\log_2\varepsilon^{-1}
\exp\{-\frac{n\varepsilon}{2}\}.
$$
Next we set $\eps_j:=2^{-j}.$
Let ${\cal J}=\{j\geq 0 : \eps_j\geq \frac{t\vee 2\log n}{n}\}$
and
$$
E:=\Bigl\{\exists j\in {\cal J}\ \exists f\in {\cal F}:\
P_n\{f\leq \delta_n^{\psi}(\eps_j)\}\leq
\eps_j\
{\rm and}\
P\{f\leq  \frac{\delta_n^{\psi}(\eps_j)}{2}\}\geq \bar A
\eps_j
\Bigr\}.
$$
We have
\begin{eqnarray}
&&
{\mathbb P}(E)
\leq B \sum_{j\in {\cal J}}
\log_2\log_2\eps_j^{-1}\exp\Bigl\{-\frac{n\eps_j}{2}\Bigr\}
\leq B\log_2\log_2\frac{n}{t\vee 2\log n}
\sum_{j\geq 0}
\exp\Bigl\{-\bigl(\frac{t}{2}\vee \log n\bigr)2^{j}\Bigr\}\leq
\nonumber
\\
&&
\leq B^{\prime}
\log_2\log_2\frac{n}{t\vee 2\log n}
\exp\Bigl\{-\bigl(\frac{t}{2}\vee \log n\bigr)\Bigr\}.
\label{bound}
\end{eqnarray}
Suppose that for some $j$ and for some $f\in {\cal F},$
$\hat \eps_n^{\psi}(t;f)\in (\eps_{j+1},\eps_j].$
On the event $E^c,$ the inequality
$P_n\{f\leq \delta_n^{\psi}(\eps_j)\}\leq \eps_j$
implies that
$P\{f\leq \delta_n^{\psi}(\eps_j)/2\}\leq \bar A \eps_j.$
Since
$$
\frac{\delta_n^{\psi}(\eps_j)}{2}=
\frac{\varphi^{-1}(\sqrt{\eps_j n})}{2\sqrt{\eps_j}}\geq
\frac{\varphi^{-1}(\sqrt{4\eps_j n})}{\sqrt{4\eps_j}}=
\delta_n^{\psi}(4\eps_j),
$$
we also have
$P\{f\leq  \delta_n^{\psi}(4\eps_j)\}\leq \bar A \eps_j,$
which implies
$P\{f\leq  \delta_n^{\psi}(8\hat \eps_n^{\psi}(f;t))\}\leq
2\bar A\hat \eps_n^{\psi}(f;t).$
Therefore, on the event $E^c,$ we get for all $f\in {\cal F},$
$\eps_n^{\psi}(f;t)\leq (2\bar A\vee 8)\hat \eps_n^{\psi}(f;t).$
It follows from (\ref{bound}) that
$$
{\mathbb P}\Bigl\{\exists f\in {\cal F}:\ \eps_n^{\psi}(f;t)\geq
(2\bar A\vee 8)\hat \eps_n^{\psi}(f;t)\Bigr\}\leq
B^{\prime}
\log_2\log_2\frac{n}{t\vee 2\log n}
\exp\Bigl\{-\bigl(\frac{t}{2}\vee \log n\bigr)\Bigr\}.
$$
Quite similarly, using the second bound of Theorem 2,
one can prove that
$$
{\mathbb P}\Bigl\{\exists f\in {\cal F}:\ \hat \eps_n^{\psi}(f;t)\geq
(2\bar A\vee 8)\eps_n^{\psi}(f;t)\Bigr\}\leq
B^{\prime}
\log_2\log_2\frac{n}{t\vee 2\log n}
\exp\Bigl\{-\bigl(\frac{t}{2}\vee \log n\bigr)\Bigr\},
$$
which implies the inequality of Theorem 1.
\qed

{\bf Proof of Theorem 2}. We follow the proof of Theorem 6
in \cite{KP}.
Define
$$
r_0:=1,\
r_{k+1}=
C\sqrt{r_k \varepsilon}\bigwedge 1
$$
where $C=c(1+{\mathbb E}D_n)$
with a sufficiently large constant $c>1$ (which will
be chosen later).
A simple induction shows that either
$C\sqrt{\eps}\geq 1$ and $r_k\equiv 1,$ or
$C\sqrt{\eps}< 1,$ and in the last case
$$
r_k=C^{1+2^{-1}+\dots +2^{-(k-1)}}\eps^{2^{-1}+\dots +2^{-k}}=
C^{2(1-2^{-k})}\eps^{1-2^{-k}}= (C\sqrt{\eps})^{2(1-2^{-k})}.
$$
Let
$\gamma_k := (\eps/r_k)^{1/2}=C^{2^{-k}-1}\eps^{2^{-k-1}}.$
Then
\begin{eqnarray}
&&
\gamma_{k}+\gamma_{k-2}+\dots +\gamma_0=
C^{-1}\bigl[C\sqrt{\eps}+(C\sqrt{\eps})^{2^{-1}}+\dots
+(C\sqrt{\eps})^{2^{-k}}\bigr]
\nonumber
\\
&&
\leq
C^{-1}(C\sqrt{\eps})^{2^{-k}}
(1-(C\sqrt{\eps})^{2^{-k}})^{-1}\leq 1/2
\label{gaga}
\end{eqnarray}
for $\eps\leq C^{-4},$ $C>2(2^{1/4}-1)^{-1}$ and
$k\leq \log_2\log_2 \eps^{-1}$
(note that $\eps\leq C^{-4}$ implies $C\sqrt{\eps}<1$).
In what follows, we fix $\eps>0$ and use only the
values of $k$ such that $k\leq \log_2\log_2 \eps^{-1}.$
%For small enough $\eps$
%(note that our choice of
%$\eps\leq C^{-4}$ implies $C\sqrt{\eps}<1$), we have
%\begin{equation}
%\gamma_0+\dots +\gamma_{k} \leq \frac{1}{2},\ k\geq 1.
%\end{equation}
Let $\delta > 0.$ Define
$$
\delta_0 = \delta,\,
\delta_k := \delta (1-\gamma_0-\dots \gamma_{k-1}),\,
\delta_{k,\frac{1}{2}}=\frac{1}{2}(\delta_k + \delta_{k+1}),\,
\ k\geq 1.
$$
Next we set ${\cal F}_0:={\cal F},$ and define recursively
$$
{\cal F}_{k+1}:=\Bigl\{f\in {\cal F}_k:
P\{f\leq \delta_{k,\frac{1}{2}}\}\leq {r_{k+1}/2}\Bigr\}.
$$
For $k\geq 0,$
define a continuous function $\varphi_k$
from ${\mathbb R}$ into $[0,1]$
such that $\varphi_k (u)=1$ for $u\leq \delta_{k,\frac{1}{2}},$
$\varphi_k (u)=0$ for $u\geq \delta_k,$  and $\varphi_k$ is linear
for $\delta_{k,\frac{1}{2}}\leq u \leq \delta_k.$
Also, for $k\geq 1,$
let $\bar \varphi_k$ be a continuous function from ${\mathbb R}$ into $[0,1]$
such that $\bar \varphi_k (u)=1$ for $u\leq \delta_k,$
$\bar \varphi_k(u)=0$ for $u\geq \delta_{k-1,\frac{1}{2}},$
and $\bar \varphi_k$ is linear for $\delta_k\leq u \leq
\delta_{k-1,\frac{1}{2}}.$
It follows from (\ref{gaga}) that
$\delta_k \in (\delta/2, \delta)$ for all
$k$ such that $1\leq k\leq \log_2\log_2 \eps^{-1}.$
Let us introduce the following function classes:
$$
{\cal G}_k :=
\bigl\{\varphi_k\circ f:f\in {\cal F}_k\bigr\},
\,\,\,k\geq 0
$$
and
$$
\bar {\cal G}_k :=
\bigl\{\bar \varphi_k\circ f:f\in {\cal F}_{k}\bigr\},
\,\,\,k\geq 1 .
$$
It follows from the definitions that, for $k\geq 1,$
$$
\sup_{g\in {\cal G}_k} Pg^2 \leq
\sup_{f\in {\cal F}_k}P\{f\leq \delta_k\}\leq
\sup_{f\in {\cal F}_k}P\{f\leq \delta_{k-1,\frac{1}{2}}\}
\leq r_k/2 \leq r_k
$$
and
$$
\sup_{g\in \bar {\cal G}_k} Pg^2 \leq
\sup_{f\in {\cal F}_{k}}P\{f\leq \delta_{k-1,\frac{1}{2}}\}
\leq r_{k}/2\leq r_k .
$$
(For $k=0,$ the first
inequality also holds since $r_0=1$).
Consider the events
$$
E^{(k)} := \Bigl\{\|P_n-P\|_{{\cal G}_{k-1}}\leq
K_1 {\mathbb E}\|P_n-P\|_{{\cal G}_{k-1}}+ K_2\sqrt{r_{k-1}\eps}+K_3\eps\Bigr\}
\bigcap
$$
$$
\bigcap \Bigl\{\|P_n-P\|_{\bar {\cal G}_k}\leq
K_1 {\mathbb E}\|P_n-P\|_{\bar{\cal G}_k}+ K_2\sqrt{r_{k}\eps}+
K_3\eps\Bigr\},\,\,k\geq 1,
$$
By concentration inequalities
of Talagrand \cite{Talagrand2,Talagrand1}
(see also \cite{Massart}), for some values of the numerical constants
$K_1, K_2, K_3>0,$
$$
{\mathbb P}((E^{(k)})^c)\leq 2e^{-\frac{n\eps}{2}}.
$$
We set $E_0=\Omega,$
$$
E_N:=\bigcap_{k=1}^N E^{(k)},
\,\,\,N\geq 1.
$$
%and
Clearly,
\begin{equation}
{\mathbb P}(E_N^c)\leq 2N e^{-\frac{n\eps}{2}}.
\label{bp}
\end{equation}
Assume, without loss of generality,
that $\eps <(2+C)^{-2},$ which implies $r_{k+1}<r_k$
and $\delta_k\in (\delta/2,\delta],\ k\geq 0.$
[If $\eps\geq (2+C)^{-2},$ the bounds of the theorem
hold with any constant $A>2+C.$]
%We will also assume that $N$ is such that
%\begin{equation}
%r_N\geq
%\varepsilon (\log_2\log_2\varepsilon^{-1})^3.
%\label{mystar}
%\end{equation}
The rest of the proof is based on the following lemma.

\begin{lemma}
Let
$$
{\cal J} :=
\Bigl\{\inf_{f\in {\cal F}}P_n\{f\leq \delta\}\leq \eps\Bigr\}.
$$
For any $N$ such that
\begin{equation}
N\leq \log_2\log_2 \eps^{-1} \mbox{ and }
r_N\geq \eps,
\label{mystars}
\end{equation}
%and (\ref{mystar}) holds.
we have on the event $E_N\bigcap {\cal J}:$
$$
(i)\ \forall f\in {\cal F}\ P_n\{f\leq \delta\}\leq \eps \Longrightarrow
f\in {\cal F}_N
$$
and
$$
(ii)\ \sup_{f\in {\cal F}_k}P_n\{f\leq \delta_k\}\leq r_k,\
0\leq k\leq N.
$$
\end{lemma}

\mdsk

{\bf Proof}. We will prove the lemma by induction
with respect to $N.$
For $N=0,$ the statement is obvious. Suppose it holds for
some $N\geq 0,$ such that $N+1$ still satisfies condition
(\ref{mystars}).
Then, on the event $E_{N}\bigcap {\cal J},$
$$
\sup_{f\in {\cal F}_k}P_n\{f\leq \delta_k\}\leq r_k,\
0\leq k\leq N
$$
and
$$
\forall f\in {\cal F}\ P_n\{f\leq \delta\}\leq \eps \Longrightarrow
f\in {\cal F}_N.
$$
Suppose that $f\in {\cal F}$ is such that
$P_n\{f\leq \delta\}\leq \eps. $ By the induction assumptions,
$f\in {\cal F}_N$ on the event $E_N.$
Hence, on the event $E_{N+1},$
\begin{eqnarray}
&&
P\{f\leq \delta_{N,\frac{1}{2}}\}
\leq P_n\{f\leq \delta_N\}+\|P_n-P\|_{{\cal G}_N}\leq
\nonumber
\\
&&
\leq \eps
+K_1 {\mathbb E}\|P_n-P\|_{{\cal G}_{N}}+ K_2\sqrt{r_{N}\eps}+
K_3\eps  .\label{e3.1}
\end{eqnarray}
Given a class ${\cal G},$ let
$$
\hat R_n ({\cal G}):=\|n^{-1}\sum_{i=1}^n \eps_i \delta_{X_i}\|_{\cal G},
$$
where $\{\eps_i\}$ is a sequence
of i.i.d. Rademacher random variables.\footnote{The random variable $\hat R_n({\cal G})$
is called \it the Rademacher complexity \rm of the class
${\cal G}.$ It was used by Koltchinskii \cite{Kolt99},
Bartlett, Boucheron and Lugosi \cite{BBL2000}, Koltchinskii and
Panchenko \cite{KP00} as a randomized complexity penalty
in learning problems}
The symmetrization inequality yields
\begin{equation}
{\mathbb E}\|P_n-P\|_{{\cal G}_{N}}\leq
2{\mathbb E}I_{E_N} {\mathbb E}_{\eps} \hat R_n({\cal G}_{N})+
2{\mathbb E}I_{E_N^c}{\mathbb E}_{\eps}\hat R_n({\cal G}_{N}).
\label{e3.2}
\end{equation}
Using the entropy inequalities for subgaussian processes
(see \cite{VandW}, Corollary 2.2.8), we get
\begin{equation}
{\mathbb E}_{\eps} \hat R_n({\cal G}_{N})\leq
\inf_{g\in {\cal G}_{N}}
{\mathbb E}_{\eps}\bigl|n^{-1}\sum_{j=1}^n \eps_j g(X_j)\bigr|
+\frac{{\rm const}}{\sqrt{n}}
\int_0^{(2\sup_{g\in {\cal G}_{N}}P_n g^2)^{1/2}}
H_{d_{P_n,2}}^{1/2}({\cal G}_{N};u)du.
\label{e3.3}
\end{equation}

{\bf Remark.} Here and in what follows in the proof
``const'' denotes a constant; its values can be different in different
places.

The induction assumption implies that on the event $E_N\bigcap {\cal J}$
$$
\inf_{g\in {\cal G}_{N}}
{\mathbb E}_{\eps}\bigl|n^{-1}\sum_{j=1}^n \eps_j g(X_j)\bigr|\leq
\inf_{g\in {\cal G}_{N}}
{\mathbb E}_{\eps}^{1/2}\bigl|n^{-1}\sum_{j=1}^n \eps_j g(X_j)\bigr|^2\leq
\frac{1}{\sqrt{n}}\inf_{g\in {\cal G}_{N}}\sqrt{P_n g^2}\leq
$$
$$
\leq \frac{1}{\sqrt{n}}
\inf_{f\in {\cal F}_{N}}\sqrt{P_n \{f\leq \delta_N\}}\leq
\frac{1}{\sqrt{n}}
\inf_{f\in {\cal F}_{N}}\sqrt{P_n \{f\leq \delta\}}\leq
\sqrt{\frac{\eps}{n}}\leq \eps,
$$
since $\eps > n^{-1}.$
Also, on the same event
%$E_N\bigcap {\cal J},$
$$
\sup_{g\in {\cal G}_{N}}P_n g^2\leq \sup_{f\in {\cal F}_N}
P_n\{f\leq \delta_N\}\leq r_N.
$$
The Lipschitz constants of $\varphi_{k-1}$ and $\bar \varphi_k$ are
bounded by
\begin{eqnarray*}
L=2(\delta_{k-1}-\delta_k)^{-1}=
2\delta^{-1}\gamma_{k-1}^{-1}=
\frac{2}{\delta}\sqrt{\frac{r_{k-1}}{\eps}},
\end{eqnarray*}
which yields
$$
d_{P_n,2}\Bigl(\varphi_N \circ f;\varphi_N\circ g\Bigr)
=
\Bigl(n^{-1}\sum_{j=1}^n \Bigl|\varphi_N (f(X_j))-
\varphi_N (g(X_j))\Bigr|^2\Bigr)^{1/2}\leq
\frac{2}{\delta}\sqrt{\frac{r_{N}}{\eps}} d_{P_n,2} (f,g).
$$
Note that for
$\eps\geq\eps_n^{\psi}(\delta)$
the inequality
$\psi(\delta \sqrt{\eps}/2)/
(\delta\sqrt{n})\leq\eps$
holds.
It follows that, on the event $E_N\bigcap {\cal J},$
\begin{eqnarray}
&&
\frac{1}{\sqrt{n}}
\int_0^{(2\sup_{g\in {\cal G}_{N}}P_n g^2)^{1/2}}
H_{d_{P_n,2}}^{1/2}({\cal G}_{N};u)du
\leq
\frac{1}{\sqrt{n}}
\int_0^{(2r_N)^{1/2}}
H_{d_{P_n,2}}^{1/2}({\cal F};
\frac{\delta \sqrt{\eps} u}{2\sqrt{r_N}})du
\nonumber
\\
&&
\leq
\frac{1}{\sqrt{n}}\frac{2\sqrt{r_N}}{\delta\sqrt{\eps}}
\int\limits_{0}^{\delta\sqrt{\eps}/2}
H_{d_{P_n,2}}^{1/2}({\cal F};v)dv
\leq
\frac{1}{\sqrt{n}}\frac{2\sqrt{r_N}}{\delta\sqrt{\eps}}
D_n\psi\bigl(\frac{\delta \sqrt{\eps}}{2}\bigr)\leq
\nonumber
\\
&&
\frac{2D_n\sqrt{r_N}}{\sqrt{\eps}}\eps=2D_n\sqrt{r_N \eps},
\label{e3.6}
\end{eqnarray}
Now (\ref{e3.3}) and (\ref{e3.6}) imply that on the same event
%$E_{N+1} \bigcap {\cal J}$
\begin{equation}
{\mathbb E}_{\eps} \hat R_n({\cal G}_{N})\leq
{\rm const}(1+D_n)\sqrt{r_N \eps}.\label{e3.7}
\end{equation}
Since
${\mathbb E}_{\eps} \hat R_n({\cal G}_{N+1})\leq 1,$
we conclude from (\ref{bp}), (\ref{e3.2}) and (\ref{e3.7})
that
$$
{\mathbb E}\|P_n-P\|_{{\cal G}_{N}}\leq
{\rm const}(1+{\mathbb E}D_n)\sqrt{r_N \eps}+
2{\mathbb P}(E_N^c)
\leq
{\rm const}(1+{\mathbb E}D_n)\sqrt{r_N \eps}+
4Ne^{-n\eps/2}.
$$
By condition (\ref{mystars}) and the fact that
$\eps\geq 2\log n /n,$ we have
$4Ne^{-n\eps/2}\leq \eps .$
Therefore,
$$
{\mathbb E}\|P_n-P\|_{{\cal G}_{N}}\leq
{\rm const}(1+{\mathbb E}D_n) \sqrt{r_N \eps}.
$$
By (\ref{e3.1}), on the event
$E_{N+1}\bigcap {\cal J}$
\begin{equation}
P\{f\leq \delta_{N,\frac{1}{2}}\}\leq
{\rm const}(1+{\mathbb E}D_n) \bigl(\eps+
\sqrt{r_N \eps}\bigr).
\end{equation}
Choosing a constant $c>0$ in the recurrent relationship defining
the sequence $\{r_k\}$ properly, we ensure that on the event
$E_{N+1}\bigcap {\cal J}$
$$
P\{f\leq \delta_{N,\frac{1}{2}}\}\leq
\frac{1}{2}C\sqrt{r_N \eps}=r_{N+1}/2.
$$
This implies that $f\in {\cal F}_{N+1}$
and the induction step for (i) is proved.

To prove (ii), note that on the event $E_{N+1}$
\begin{eqnarray}
&&
\sup_{f\in {\cal F}_{N+1}} P_n\{f\leq \delta_{N+1}\}\leq
\sup_{f\in {\cal F}_{N+1}} P\{f\leq \delta_{N,\frac{1}{2}}\}+
\|P_n-P\|_{\bar {\cal G}_{N+1}}\leq
\nonumber
\\
&&
\leq r_{N+1}/2 +
K_1 {\mathbb E}\|P_n-P\|_{\bar {\cal G}_{N+1}}+ K_2\sqrt{r_{N+1}\eps}+
K_3\eps .\label{ee3.1}
\end{eqnarray}
Using the symmetrization inequality, we get
\begin{equation}
{\mathbb E}\|P_n-P\|_{\bar {\cal G}_{N+1}}\leq
2{\mathbb E}I_{E_N} {\mathbb E}_{\eps} \hat R_n(\bar {\cal G}_{N+1})+
2{\mathbb E}I_{E_N^c}{\mathbb E}_{\eps}\hat R_n(\bar {\cal G}_{N+1}).
\label{ee3.2}
\end{equation}
Similarly to (\ref{e3.3})
\begin{equation}
{\mathbb E}_{\eps} R_n(\bar {\cal G}_{N+1})\leq
\inf_{g\in \bar {\cal G}_{N+1}}
{\mathbb E}_{\eps}\bigl|n^{-1}\sum_{j=1}^n \eps_j g(X_j)\bigr|
+\frac{{\rm const}}{\sqrt{n}}
\int_0^{(2\sup_{g\in \bar {\cal G}_{N+1}}P_n g^2)^{1/2}}
H_{d_{P_n,2}}^{1/2}(\bar {\cal G}_{N+1};u)du.
\label{ee3.3}
\end{equation}
It follows from (i) that on the event $E_{N+1}\bigcap {\cal J}$
$$
\inf_{g\in \bar {\cal G}_{N+1}}
{\mathbb E}_{\eps}\bigl|n^{-1}\sum_{j=1}^n \eps_j g(X_j)\bigr|\leq
\inf_{g\in \bar {\cal G}_{N+1}}
{\mathbb E}_{\eps}^{1/2}\bigl|n^{-1}\sum_{j=1}^n \eps_j g(X_j)\bigr|^2\leq
\frac{1}{\sqrt{n}}\inf_{g\in \bar {\cal G}_{N+1}}\sqrt{P_n g^2}\leq
$$
$$
\leq \frac{1}{\sqrt{n}}
\inf_{f\in {\cal F}_{N+1}}\sqrt{P_n \{f\leq \delta_{N,\frac{1}{2}}\}}\leq
\frac{1}{\sqrt{n}}
\inf_{f\in {\cal F}_{N+1}}\sqrt{P_n \{f\leq \delta\}}\leq
\sqrt{\frac{\eps}{n}}\leq \eps.
$$
The induction assumption implies that on the event
$E_{N+1}\bigcap {\cal J}$
$$
\sup_{g\in \bar {\cal G}_{N+1}}P_n g^2\leq \sup_{f\in {\cal F}_N}
P_n\{f\leq \delta_{N,\frac{1}{2}}\}\leq r_N.
$$
Since the Lipschitz constant of $\bar \varphi_k$ is
bounded by $\frac{2}{\delta}\sqrt{\frac{r_{k-1}}{\eps}},$
we have
$$
d_{P_n,2}\Bigl(\bar \varphi_{N+1}\circ f;
\bar \varphi_{N+1}\circ g\Bigr)
=
\Bigl(n^{-1}\sum_{j=1}^n \Bigl|\bar \varphi_{N+1} \circ f(X_j)-
\bar \varphi_{N+1} \circ g(X_j) \Bigr|^2\Bigr)^{1/2}\leq
\frac{2}{\delta}\sqrt{\frac{r_{N}}{\eps}} d_{P_n,2}(f,g).
$$
Similarly to (\ref{e3.6}),
we have on the event $E_{N+1}\bigcap {\cal J},$
\begin{eqnarray}
&&
\frac{1}{\sqrt{n}}
\int_0^{(2\sup_{g\in \bar {\cal G}_{N+1}}P_n g^2)^{1/2}}
H_{d_{P_n,2}}^{1/2}(\bar {\cal G}_{N+1};u)du
\leq
\frac{1}{\sqrt{n}}
\int_0^{(2r_N)^{1/2}}
H_{d_{P_n,2}}^{1/2}({\cal F};
\frac{\delta \sqrt{\eps} u}{2\sqrt{r_N}})du
\nonumber
\\
&&
\leq
\frac{1}{\sqrt{n}}\frac{2\sqrt{r_N}}{\delta\sqrt{\eps}}
\int\limits_{0}^{\delta\sqrt{\eps}/2}
H_{d_{P_n,2}}^{1/2}({\cal F};v)dv
\leq
\frac{1}{\sqrt{n}}\frac{2\sqrt{r_N}}{\delta\sqrt{\eps}}
D_n\psi\bigl(\frac{\delta \sqrt{\eps}}{2}\bigr)\leq
\nonumber
\\
&&
\frac{2D_n\sqrt{r_N}}{\sqrt{\eps}}\eps=2D_n\sqrt{r_N \eps}.
\label{ee3.6}
\end{eqnarray}
Combining all the bounds, we prove that on the same event
%$E_{N+1}\bigcap {\cal J}$
\begin{equation}
\sup_{f\in {\cal F}_{N+1}}P_n\{f\leq \delta_{N+1}\}\leq
\frac{r_{N+1}}{2}+
{\rm const}(1+{\mathbb E}D_n) \sqrt{r_N \eps}.
\end{equation}
Choosing a constant $c>0$ in the recurrent relationship defining
the sequence $\{r_k\}$ properly, we get on the event
$E_{N+1}\bigcap {\cal J}$
$$
\sup_{f\in {\cal F}_{N+1}}P_n\{f\leq \delta_{N+1}\}\leq
C\sqrt{r_N \eps}=r_{N+1},
$$
which completes the proof of (ii)
and of the lemma.
\qed

To complete the proof of the theorem, note that
the choice of
$N=[\log_2\log_2 \eps^{-1}]$
implies that $r_{N+1}\leq c\eps$ for some $c>0.$
Indeed, if we introduce $s_k=r_k/C$ and $\eps_1=C\eps$ then
$s_{k+1}=\sqrt{s_k\eps}$ and $s_0=C^{-1}\leq 1.$
It is easy to see that
$s_N\leq \eps_1^{1-2^{-N}}\leq 2\eps_1$
for $N\geq \log_2\log_2 \eps_1^{-1},$
and, hence,
$r_N\leq C^2\eps=\bar A\eps.$

The proof of the second inequality
is similar with minor modifications.
\qed

To prove Theorem 3, we need the following statement, which
seems to be well known, but we have not found the precise reference
and give the proof here for completeness.

Let
$$
{\rm conv}_d({\cal H}):= \Bigl\{\sum_{j=1}^d \lambda_j h_j:
\lambda_j\in {\mathbb R},
\sum_{j=1}^d|\lambda_j|\leq 1,\ h_j\in {\cal H}\Bigr\}.
$$

\begin{lemma} Let ${\cal H}$ be a class of functions from
$(S,{\cal A})$ into ${\mathbb R}.$ Let $Q$ be a probability measure on $(S,{\cal A})$
such that
$$
\bar H:=\sup_{h\in {\cal H}}(Qh^2)^{1/2}<+\infty .
$$
The following bound holds for all $d\geq 1$ and $\eps>0:$
$$
N_{d_{Q,2}}\Bigl({\rm conv}_d({\cal H}), (1+\bar H)\eps\Bigr)
\leq \left(\frac{2e^2 N_{d_{Q,2}}({\cal H}, \eps )(d^{\prime}+4\eps^{-2})}
{{d^{\prime}}^{2}}\right)^{d^{\prime}},
$$
where $d^{\prime}=d\wedge N_{d_{Q,2}}({\cal H}, \eps ).$
\end{lemma}

{\bf Proof}. First note that if
${\cal H}^{\prime}:={\cal H}\bigcup \bigl\{h:-h\in {\cal H}\bigr\},$
then ${\rm conv}_d({\cal H}^{\prime})={\rm conv}_d({\cal H})$
and
$$
N_{d_{Q,2}}({\cal H}^{\prime};\eps)\leq
2 N_{d_{Q,2}}({\cal H};\eps).
$$
Thus, it's enough to show that for a class ${\cal H},$ such that
$h\in {\cal H}$ implies $-h\in {\cal H},$ we have
$$
N_{d_{Q,2}}\Bigl({\rm conv}_d({\cal H}), (1+\bar H)\eps\Bigr)
\leq \left(\frac{e^2 N_{d_{Q,2}}({\cal H}, \eps )(d+4\eps^{-2})}
{d^2}\right)^{d}.
$$
For such a class we have
$$
{\rm conv}_d({\cal H}):= \Bigl\{\sum_{j=1}^d \lambda_j h_j:
\lambda_j\geq 0,
\sum_{j=1}^d \lambda_j\leq 1,\ h_j\in {\cal H}\Bigr\}.
$$
Note that if $\sum_{j}|\lambda_j|\leq 1,$ then
$$
d_{Q,2}\Bigl(\sum_j \lambda_j h_j;\sum_j \lambda_j h_j^{\prime}\Bigr)=
\Bigl\|\sum_j \lambda_j (h_j-h_j^{\prime})\Bigr\|_{L_2(Q)}\leq
$$
$$
\leq \sum_{j}|\lambda_j|\max_{j}\bigl\|h_j-h_j^{\prime}\bigr\|_{L_2(Q)}\leq
\max_{j}\bigl\|h_j-h_j^{\prime}\bigr\|_{L_2(Q)}.
$$
It follows that if ${\cal H}_{\eps}$ is an $\eps$-net of ${\cal H},$
then a $\delta$-net of ${\rm conv}_d({\cal H}_{\eps})$ is an $\eps+\delta$-net
of ${\rm conv}_d({\cal H}).$ This observation allows us to reduce
the proof of the lemma to the case when ${\cal H}$ is a finite
class. In this case we want to show that
$$
N_{d_{Q,2}}\Bigl({\rm conv}_d({\cal H}), \bar H\eps\Bigr)
\leq \left(\frac{e^2 {\rm card}({\cal H})(d+4\eps^{-2})}
{d^2}\right)^{d}.
$$
To this end, we use the idea of
B. Maurey, see \cite{Pisier,VandW}.
Let $N:={\rm card}({\cal H}).$
Consider some representation of a function
$
f=\sum_{i=1}^{N} \lambda_i h_i \in {\rm conv}_d({\cal H}).
$
We assume that $\lambda_j\geq 0,$ $\sum_{j}\lambda_j\leq 1,$
and at most $d^{\prime}$ of the coefficients are not equal to $0.$
Consider an i.i.d. sequence of random variables $Y_j,\,j=1,\ldots,k$
taking values in ${\cal H}\cup \{0\}$
such that $P(Y_j=h_i)=\lambda_i$ for $i=1,\dots, N$ and
$P(Y_j=0)=1-\sum_{i=1}^N\lambda_i.$
(We simply add the probabilities when the same function $h$
corresponds to several weights $\lambda_i$ with different
indices). We have
\begin{eqnarray*}
&&
{\mathbb E}\|k^{-1}\sum_{j=1}^k Y_j -
\sum_{i=1}^{N}\lambda_i h_i\|_{Q,2}^2 =
{\mathbb E}\|k^{-1}\sum_{j=1}^k Y_j - {\mathbb E}Y_1\|_{Q,2}^2 \leq
\\
&&
\leq
\frac{1}{k}{\mathbb E}\|Y_1 - {\mathbb E}Y_1\|_{Q,2}^2
\leq 4{\bar H}^2 k^{-1}.
\end{eqnarray*}
If we set $k=4\eps^{-2},$ then with probability $1$ there exists
a realization $\bar{Y}_k = k^{-1}\sum_{j=1}^k Y_j$
such that
$$
\|\bar{Y}_k - \sum_{i=1}^{N}\lambda_i h_i\|_{Q,2}
\leq \eps \bar H.
$$
In order to compute the bound for
the $\bar H\eps-$covering
number we have to calculate the number of possible
realizations of $k^{-1}\sum_{j=1}^k Y_j.$
A simple combinatorics shows that this number does not exceed
${N\choose {d^{\prime}}}{{d^{\prime}+k} \choose k}.$
Next we use the following bound, which holds for all
$1\leq d\leq N:$
$$
{N\choose d}{d+k \choose k}\leq \left(
\frac{e^2N(d+k)}{d^2}
\right)^{d}.
$$
To prove the bound, first assume that $d<N.$
Then one can check using
Stirling's formula that
\begin{eqnarray*}
&&
\frac{N!}{d!(N-d)!}
\frac{(d+k)!}{d!k!}\leq
\frac{n^n}{d^d (N-d)^{N-d}}
\frac{(d+k)^{d+k}}{k^k d^d}
\\
&&
\leq
\left(\frac{N(d+k)}{d^2}\right)^d
\left(1+\frac{d}{N-d}\right)^{N-d}
\left(1+\frac{d}{k}\right)^{k}\leq
\left(\frac{e^2 N(d+k)}{d^2}\right)^d.
\end{eqnarray*}
The case when $d=N$ can be considered similarly.
The bound immediately implies the result.
\qed

{\bf Proof of Theorem 3.}
Let us fix $\delta\in (0,1/2].$
For any function $f$ we denote $d(f):=d(f,\bar \Delta),$
where $\bar \Delta$ is such that
the infimum in the definition (\ref{efd})
is attained at $\bar \Delta.$
For a fixed $\delta$
we consider a partition of ${\cal F}$ into two classes
${\cal F}_1^{\delta}$ and
${\cal F}_2^{\delta} = {\cal F}\setminus {\cal F}_1^{\delta},$
where ${\cal F}_1^{\delta}:=\{f : d(f)=0\}$
(note that $d(f)$ depends on $\delta$).
In the first four steps of the proof we will
deal with ${\cal F}_2^{\delta}$
and we will assume only that the class $\cal H$
has a square integrable envelope $H.$

{\it Step 1}.
Let $1\leq d\leq n.$
Denote
$$
\eps_n(d;\delta ;\Delta):=
\Bigr[\frac{d}{n}\Bigl(\log\frac{1}{\delta}+
\log{\frac{ne^2}{d}}\Bigr)+
\Bigl(\frac{\Delta}{\delta}\Bigr)^{\frac{2\alpha}{\alpha+2}}
n^{-\frac{2}{\alpha+2}}\Bigr]\bigvee \frac{2\log n}{n}.
$$
Let
${\cal F}_{d,\Delta}:=\{f\in {\cal F}_2^{\delta}: d(f;\Delta)\leq d\}.$
We start by proving (with some constants $A,B>0$)
the following inequality:
\begin{eqnarray}
&&
{\mathbb P}\Bigl\{\exists f\in {\cal F}_{d,\Delta}\
P_n\{f\leq \delta\}\leq
\eps_n(d;\delta;\Delta)\
{ \rm and }\
P\{f\leq  \frac{\delta}{2}\}\geq A
\eps_n(d;\delta;\Delta)
\Bigr\}\leq
\nonumber
\\
&&
\leq B
\Bigl(\frac{\delta d}{n}\Bigr)^{d/4}
\exp\Bigl\{-\frac{1}{4}\bigl(\sqrt{n}
\frac{\Delta}{\delta}\bigr)^{2\alpha/(\alpha+2)}\Bigr\}.
\label{st1}
\end{eqnarray}
Clearly, we can and do assume that $\eps_n(d;\delta;\Delta)\leq 1.$
To prove (\ref{st1}), we bound the random entropy
$H_{d_{P_n,2}}({\cal F}_{d,\Delta};\eps)$ of the
class ${\cal F}_{d,\Delta}$ the following way:
\begin{equation}
H_{d_{P_n,2}}({\cal F}_{d,\Delta};\eps)\leq
K(1+P_n H^2)\Bigl[d \log\frac{e}{\eps}
+\Bigl(\frac{\Delta}{\eps}\Bigr)^{\alpha}\Bigr]
\ {\rm for}\ \eps \leq 1
\label{ent1}
\end{equation}
with some constant $K>0.$
The last bound follows from the observation that each function
$f\in {\cal F}_{d,\Delta}$ can be represented as $f=f_1+f_2,$
where
$$
f_1\in {\cal F}_d :=
{\rm conv}_d({\cal H})= \Bigl\{\sum_{j=1}^d \lambda_j h_j:
\lambda_j\in {\mathbb R},
\sum_{j=1}^d|\lambda_j|\leq 1,\ h_j\in {\cal H}\Bigr\}
$$
and
$$
f_2\in {\cal F}_{\Delta}:=\Delta\ {\rm conv}({\cal H}).
$$
Hence, by simple combining of $\eps$-coverings for the classes
${\cal F}_d$ and ${\cal F}_{\Delta},$ we get
$$
H_{d_{P_n,2}}({\cal F}_{d,\Delta};\eps)\leq
H_{d_{P_n,2}}({\cal F}_{d};\eps/2)+
H_{d_{P_n,2}}({\cal F}_{\Delta};\eps/2).
$$
Then, a routine application of Lemma 2 and (\ref{envc})
implies
$$
H_{d_{P_n,2}}({\cal F}_{d};\eps/2)\leq
Kd\log{\frac{e(1+P_n H^2)}{\eps}}\ {\rm for}\ \eps \leq
2(P_n H^2)^{1/2}
$$
(note that for $\eps> 2(P_n H^2)^{1/2}$ we easily get
$H_{d_{P_n,2}}({\cal F}_{d};\eps/2)=0$).
For $\eps \leq 1$ this implies
$$
H_{d_{P_n,2}}({\cal F}_{d};\eps/2)\leq
Kd\Bigl[\log\frac{e}{\eps}+\log(1+P_n H^2)\Bigr]
\leq Kd\Bigl[\log\frac{e}{\eps}+P_n H^2\Bigr]
\leq Kd(1+P_n H^2)\log\frac{e}{\eps}.
$$
By the bound on the entropy of the symmetric convex hull
(see \cite{VandW})
$$
H_{d_{P_n,2}}({\cal F}_{\Delta};\eps/2)=
H_{d_{P_n,2}}\bigl({\cal F};{\frac{\eps}{2\Delta}}\bigr)\leq
K(1+P_n H^2)^{\alpha/4}\Bigl({\frac{\Delta}{\eps}}\Bigr)^{\alpha}\leq
K(1+P_n H^2)\Bigl({\frac{\Delta}{\eps}}\Bigr)^{\alpha},
$$
which implies (\ref{ent1}).

Next we are using margin-type bounds on generalization error under
random entropy conditions (see Section 2, Theorem 2).
Clearly, from (\ref{ent1}),
we get the following bound on Dudley's entropy integral:
$$
\int_{0}^x H_{d_{P_n,2}}^{1/2}({\cal F};\eps)d\eps \leq
K(1+P_n H^2)^{1/2} \bar \psi (x),
$$
where $\bar \psi$ is a concave nondecreasing function
such that for $x\in [0,1]$
$$
\bar \psi (x)=\Bigl(x\Bigl(d\log{\frac{e}{x}}\Bigr)^{1/2}
+\Delta^{\alpha/2}x^{1-\alpha/2}\Bigr)
$$
with some constant $K>0.$
Let
$$
\psi_1(x):=x\Bigl(d\log{\frac{e}{x}}\Bigr)^{1/2} ,\
\psi_2(x):=\Delta^{\alpha/2}x^{1-\alpha/2},\
\psi(x):=(\psi_1(x)+\psi_2(x))/2.
$$
Let us first consider the equation
$\eps = \psi_1 (\delta\sqrt{\eps})/(\delta \sqrt{n}),$
which can be written as
$\eps=\frac{d}{n}\log\frac{e}{\delta\sqrt{\eps}}.$
If $\eps=\frac{d}{n}x^2$ then
$$
xe^{x^2}=\Bigl(\frac{n}{d}\Bigr)^{1/2}\frac{e}{\delta}.
$$
For $d\leq n$ and $\delta \le 1,$ it means that $xe^{x^2}\geq 1,$
and, therefore,
$$
e^{x^2 -1}\leq \Bigl(\frac{n}{d}\Bigr)^{1/2}\frac{e}{\delta},
$$
or,
$$
\eps=\frac{d}{n}x^2
\leq \frac{d}{n}\Bigl[1+
\log\Bigl(\Bigl(\frac{n}{d}\Bigr)^{1/2}\frac{e}{\delta}\Bigr)\Bigr]
\leq
\frac{d}{n} \log\frac{ne^2}{d\delta}\leq
\eps_n(d;\delta;\Delta)\leq 1.
$$
[One can notice that in the case when $d$ becomes
significantly greater then $n,$ for example,
if $(nd^{-1})^{1/2}\delta^{-1}\leq 1$
then $x\leq 1$ and $xe^{x^2}\leq e x,$
which implies that $\eps\geq \delta^{-2}$
and the bound of the theorem becomes useless.
This explains why in the definition of $\eps_n(f;\delta)$
we minimize over $d(f,\Delta)\leq n.$]

The solution of the equation
$\eps = \psi_2 (\delta\sqrt{\eps})/(\delta \sqrt{n})$
is equal to
$$
\eps^{(2)}:=
\Bigl(\frac{\Delta}{\delta}\Bigr)^{\frac{2\alpha}{\alpha+2}}
n^{-\frac{2}{\alpha+2}}.
$$
Finally, it is easy to bound
the solution of the equation
$\eps =\psi(\delta\sqrt{\eps})/(\delta \sqrt{n})$
from above by $\eps^{(1)}+\eps^{(2)}.$
Therefore, the solution of the last equation is also bounded
from above by $\eps_n(d;\delta;\Delta).$
This allows us to use the bound of Theorem 2 to get the following
inequality:
$$
{\mathbb P}\Bigl\{\exists f\in {\cal F}_{d,\Delta}\
P_n\{f\leq \delta\}\leq
\eps_n(d;\delta;\Delta)\
{\rm and}\
P\{f\leq  \frac{\delta}{2}\}\geq A
\eps_n(d;\delta;\Delta)
\Bigr\}\leq
$$
$$
\leq B \log_2\log_2\eps_n(d;\delta;\Delta)^{-1}
\exp\{-\frac{n\eps_n(d;\delta;\Delta)}{2}\}.
$$
Since, for $\eps := \eps_n(d;\delta;\Delta),$ we have
$\eps \geq \frac{2\log n}{n},$ it follows that for $n\geq 3,$
$$
\frac{1}{\eps}\log \log_2 \log_2 \frac{1}{\eps}\leq n/4,
$$
which implies
\begin{equation}
B \log_2\log_2\eps_n(d;\delta;\Delta)^{-1}
\exp\{-\frac{n\eps_n(d;\delta;\Delta)}{2}\}\leq
B\exp\{-\frac{n\eps_n(d;\delta;\Delta)}{4}\}.
\label{ecomp}
\end{equation}
A simple computation shows that
$$
\exp\{-\frac{n\eps_n(d;\delta;\Delta)}{4}\}\leq
\Bigl(\frac{\delta d}{n}\Bigr)^{d/4}
\exp\Bigl\{-\frac{1}{4}\bigl(\sqrt{n}
\frac{\Delta}{\delta}\bigr)^{2\alpha/(\alpha+2)}\Bigr\},
$$
which implies (\ref{st1})

{\it Step 2}. Next we show that with some constants $A,B\geq 1,$
$\delta \leq 1/2$ and $\Delta\geq\delta n^{-1/2}$
\begin{eqnarray}
&&
{\mathbb P}\Bigl\{\exists f\in {\cal F}_2^{\delta}\
P_n\{f\leq \delta\}\leq
\eps_n(d(f;\Delta);\delta;\Delta)\
{ \rm and }\
P\{f\leq  \frac{\delta}{2}\}\geq A
\eps_n(d(f;\Delta);\delta;\Delta)
\Bigr\}\leq
\nonumber
\\
&&
\leq B \delta^{1/8}\Delta^{1/8}
\exp\Bigl\{-\frac{1}{4}\bigl(\sqrt{n}
\frac{\Delta}{\delta}\bigr)^{2\alpha/(\alpha+2)}\Bigr\},
\label{st2}
\end{eqnarray}
where it's understood that if $d=d(f;\Delta)>n$
then $\eps_n(d;\delta;\Delta)=1.$
Indeed, using (\ref{st1}), we have for $\delta\leq 1/2$
$$
{\mathbb P}\Bigl\{\exists f\in {\cal F}_2^{\delta}\
P_n\{f\leq \delta\}\leq
\eps_n(d(f;\Delta);\delta;\Delta)\
{\rm and}\
P\{f\leq  \frac{\delta}{2}\}\geq A
\eps_n(d(f;\Delta);\delta;\Delta)
\Bigr\}\leq
$$
$$
\leq {\mathbb P}\Bigl\{\exists d\leq n\ \exists f\in {\cal F}_2^{\delta}\
d(f;\Delta)=d, P_n\{f\leq \delta\}\leq
\eps_n(d;\delta;\Delta)\
{\rm and}\
P\{f\leq  \frac{\delta}{2}\}\geq A
\eps_n(d;\delta;\Delta)
\Bigr\}\leq
$$
$$
\leq \sum_{d=1}^{n}{\mathbb P}\Bigl\{\exists f\in {\cal F}_{d,\Delta}\
 P_n\{f\leq \delta\}\leq
\eps_n(d;\delta;\Delta)\
{ \rm and }\
P\{f\leq  \frac{\delta}{2}\}\geq A
\eps_n(d;\delta;\Delta)
\Bigr\}\leq
$$
$$
\leq
B\sum_{d=1}^n\Bigl(\frac{\delta d}{n}\Bigr)^{d/4}
\exp\Bigl\{-\frac{1}{4}\bigl(\sqrt{n}
\frac{\Delta}{\delta}\bigr)^{2\alpha/(\alpha+2)}\Bigr\}.
$$
One can easily check that for $d\leq n/(e\delta)$
(increasing $A$ we can assume that it holds)
the expression $(\delta d/n)^{d/4}$ is
decreasing in $d$ and, therefore, for any $k\leq n/e$
$$
\sum_{d=1}^n\Bigl(\frac{\delta d}{n}\Bigr)^{d/4}\leq
k\Bigl(\frac{\delta}{n}\Bigr)^{1/4}+
\sum_{d=k+1}^n\Bigl(\frac{\delta d}{n}\Bigr)^{d/4}\leq
k\Bigl(\frac{\delta}{n}\Bigr)^{1/4}+
\delta^{k/4}.
$$
Optimizing over $k$ we take $k=\log n /\log \delta^{-1} +1$
to get
$$
k\Bigl(\frac{\delta}{n}\Bigr)^{1/4}+\delta^{k/4}
\leq
2\Bigl(\frac{\log n}{\log \delta^{-1}}+1\Bigr)
\Bigl(\frac{\delta}{n}\Bigr)^{1/4}\leq
\delta^{1/8}\Delta^{1/8},
$$
where the last inequality holds under the assumption that
$\Delta\geq\delta n^{-1/2}.$

{\it Step 3}. Our next goal is to prove that with some
constants $A,B>1$ and for $0<t<n^{\alpha/(2+\alpha)}$
\begin{eqnarray}
&&
{\mathbb P}\Bigl\{\exists f\in {\cal F}_2^{\delta}\
P_n\{f\leq \delta\}\leq
\eps_n(f;\delta)\
{\rm and}\
P\{f\leq  \frac{\delta}{2}\}\geq A
\inf_{\Delta\geq \delta n^{-1/2}
t^{\frac{1}{\alpha}+\frac{1}{2}}}
\eps_n(d(f;\Delta);\delta;\Delta)
\Bigr\}\leq
\nonumber
\\
&&
\leq B \delta^{1/8} e^{-t/4}
\label{st3}
\end{eqnarray}
Let $\Delta_j:=2^{-j},\ j\geq 0.$
Let ${\cal J}=\{j:\Delta_j\geq \delta n^{-1/2}
t^{\frac{1}{\alpha}+\frac{1}{2}}\}.$
Note that the condition $t<n^{\alpha/(2+\alpha)}$
guarantees that ${\cal J}\not = \emptyset .$
Using (\ref{st2}), we get
\begin{eqnarray*}
&&
{\mathbb P}\Bigl\{\exists f\in {\cal F}_2^{\delta}\
P_n\{f\leq \delta\}\leq
\eps_n(f;\delta)\
{\rm and}\
P\{f\leq  \frac{\delta}{2}\}\geq A
\inf_{\cal J}
\eps_n(d(f;\Delta_j);\delta;\Delta_j)
\Bigr\}\leq
\\
&&
\leq {\mathbb P}\Bigl\{\exists f\in {\cal F}_2^{\delta}\
\exists j\in {\cal J}\
P_n\{f\leq \delta\}\leq
\eps_n(f;\delta)\
{ \rm and }\
P\{f\leq  \frac{\delta}{2}\}\geq A
\eps_n(d(f;\Delta_j);\delta;\Delta_j)
\Bigr\}\leq
\\
&&
\leq
\sum_{\cal J}
{\mathbb P}\Bigl\{\exists f\in {\cal F}_2^{\delta}\
P_n\{f\leq \delta\}\leq
\eps_n(f;\delta)\
{\rm and}\
P\{f\leq  \frac{\delta}{2}\}\geq A
\eps_n(d(f;\Delta_j);\delta;\Delta_j)
\Bigr\}\leq
\\
&&
\leq B
\sum_{\cal J}
\delta^{1/8} \Delta_j^{1/8}
\exp\Bigl\{-\frac{1}{4}\bigl(\sqrt{n}
\frac{\Delta_j}{\delta}\bigr)^{2\alpha/(\alpha+2)}\Bigr\}\leq
B' \delta^{1/8} e^{-t/4}.
\end{eqnarray*}
To complete the proof of (\ref{st3}), note that
for $\Delta \in (\Delta_{j+1},\Delta_j]$ we have
$$
\frac{d(f;\Delta_j)}{n}\Bigl(\log\frac{1}{\delta}+
\log{\frac{ne}{d(f;\Delta_j)}}\Bigr)\leq
\frac{d(f;\Delta)}{n}\Bigl(\log\frac{1}{\delta}+
\log{\frac{ne}{d(f;\Delta)}}\Bigr),
$$
$$
\Bigl(\frac{\Delta_j}{\delta}\Bigr)^{\frac{2\alpha}{\alpha+2}}
n^{-\frac{2}{\alpha+2}}\leq
2^{\frac{2\alpha}{(\alpha+2)}}
\Bigl(\frac{\Delta}{\delta}\Bigr)^{\frac{2\alpha}{\alpha+2}}
n^{-\frac{2}{\alpha+2}},\,\,\,
\log\log\frac{2}{\Delta_j}\leq \log\log\frac{2}{\Delta},
$$
which implies
$\eps_n(f;\Delta_j;\delta)\leq
2^{2\alpha/(\alpha+2)}\eps_n(f;\Delta;\delta)$
and, therefore,
$$
\inf_{\cal J}
\eps_n(d(f;\Delta_j);\delta;\Delta_j)\leq
2^{2\alpha/(\alpha+2)}
\inf_{\Delta\geq \delta n^{-1/2}
t^{\frac{1}{\alpha}+\frac{1}{2}}}
\eps_n(d(f;\Delta);\delta;\Delta),
$$
and (\ref{st3}) follows.

{\it Step 4}. Now we prove that
for some constants $A,B>1$ and for all $0<t<n^{\alpha/{2+\alpha}}$
\begin{eqnarray}
&&
{\mathbb P}\Bigl\{\exists f\in {\cal F}_2^{\delta}\
P_n\{f\leq \delta\}\leq
\eps_n(f;\delta)\
{\rm and}\
P\{f\leq  \frac{\delta}{2}\}\geq A
\bigl(\eps_n(f;\delta)+\frac{t}{n}\bigr)
\Bigr\}\leq
\nonumber
\\
&&
\leq B \delta^{1/8} e^{-t/4}
\label{st4}
\end{eqnarray}
Because of (\ref{st3}), it is enough to show that
\begin{equation}
\inf_{\Delta\geq \delta n^{-1/2}
t^{\frac{1}{\alpha}+\frac{1}{2}}, \Delta\in \Delta_f}
\eps_n(d(f;\Delta);\delta;\Delta)\leq
\eps_n(f;\delta)+\frac{t}{n}.
\label{st4'}
\end{equation}
Since $d(f;\Delta)$ is a decreasing function of $\Delta ,$
the set $\Delta_f$ is an interval of the form $[c,1]$
for some $c\leq 1.$ Let
$\Delta_0 := \delta n^{-1/2}t^{\frac{1}{\alpha}+\frac{1}{2}}.$
If $\Delta_0\not\in \Delta_f,$ then (\ref{st4'}) clearly holds.
Otherwise, suppose that the infimum in the definition of $\eps_n(f;\delta)$
is attained at $\Delta = \bar \Delta.$
If $\bar \Delta \geq \Delta_0,$ then (\ref{st4'}) is also obvious.
In the case when $\bar \Delta < \Delta_0,$
note that
$$
\Bigl(\frac{\Delta_0}{\delta}\Bigr)^{\frac{2\alpha}{\alpha+2}}
n^{-\frac{2}{\alpha+2}}=\frac{t}{n}
$$
and the function
$
\frac{d(f;\Delta)}{n}\Bigl(\log\frac{1}{\delta}+
\log{\frac{ne^2}{d(f;\Delta)}}\Bigr)
$
is decreasing in $\Delta.$ Therefore,
$$
\inf_{\Delta\geq \delta n^{-1/2}
t^{\frac{1}{\alpha}+\frac{1}{2}}, \Delta\in \Delta_f}
\eps_n(d(f;\Delta);\delta;\Delta)\leq
\eps_n(d(f;\Delta_0);\delta;\Delta_0)\leq
\frac{d(f;\bar \Delta)}{n}\Bigl(\log\frac{1}{\delta}+
\log{\frac{ne^2}{d(f;\bar \Delta)}}\Bigr)+\frac{t}{n}\leq
$$
$$
\leq
\eps_n(d(f;\bar \Delta);\delta;\bar \Delta)+\frac{t}{n}
\leq
\eps_n(f;\delta)+\frac{t}{n},
$$
which proves (\ref{st4'}).

{\it Step 5}. To complete the proof of the theorem,
define the following event
$$
E:=
\Bigl\{\exists f\in {\cal F}\ \exists \delta \in (0,1):
P_n\{f\leq \delta\}\leq
\eps_n(f;\delta)\
{\rm and}\
P\{f\leq  \frac{\delta}{4}\}\geq A
\bigl(\eps_n(f;\frac{\delta}{2})+\frac{t}{n}\bigr)
\Bigr\}.
$$
Obviously, $E=E_1\bigcup E_2,$ where
$$
E_1:=
\Bigl\{ \exists \delta\in (0,1)\ \exists f\in {\cal F}_1^{\delta}\ :
P_n\{f\leq \delta\}\leq
\eps_n(f;\delta)\
{\rm and}\
P\{f\leq  \frac{\delta}{4}\}\geq A
\bigl(\eps_n(f;\frac{\delta}{2})+\frac{t}{n}\bigr)
\Bigr\},
$$
$$
E_2:=
\Bigl\{ \exists \delta \in (0,1) \ \exists f\in {\cal F}_2^{\delta}\ :
P_n\{f\leq \delta\}\leq
\eps_n(f;\delta)\
{\rm and}\
P\{f\leq  \frac{\delta}{4}\}\geq A
\bigl(\eps_n(f;\frac{\delta}{2})+\frac{t}{n}\bigr)
\Bigr\}.
$$

We set $\delta_j:=2^{-j},\,j\geq 0$ and
$$
\bar E_2:=
\Bigl\{ \exists j\geq 0 \ \exists f\in {\cal F}_2^{\delta_j}\ :
P_n\{f\leq \delta_j\}\leq
\eps_n(f;\delta_j)\
{\rm and}\
P\{f\leq  \frac{\delta_j}{2}\}\geq A
\bigl(\eps_n(f;\delta_j)+\frac{t}{n}\bigr)
\Bigr\}.
$$
It is easily seen that $E_2\subset \bar E_2.$
It follows from (\ref{st4}) that
\begin{eqnarray*}
&&
{\mathbb P}(E_2)\leq {\mathbb P}(\bar E_2)\leq
\sum_{j=0}^{\infty}{\mathbb P}
\Bigl\{\exists f\in {\cal F}_2^{\delta_j}:\
P_n\{f\leq \delta_j\}\leq
\eps_n(f;\delta_j)\
\\
&&
{\rm and}\
P\{f\leq  \frac{\delta_j}{2}\}\geq A
\bigl(\eps_n(f;\delta_j)+\frac{t}{n}\bigr)
\Bigr\}\leq \sum_{j=0}^{\infty}B
\delta_j^{1/8}e^{-t/4}\leq B' e^{-t/4}.
\end{eqnarray*}
If $f=\sum \lambda_i h_i \in {\cal F}_1^{\delta}$
for some $\delta$ then
$$
\eps_n(f,\delta)=
\Bigl(\frac{\Delta(f)}{\delta}\Bigr)^{\frac{2\alpha}{2+\alpha}}
n^{-\frac{2}{2+\alpha}}\bigvee \frac{2\log n}{n}.
$$
where $\Delta(f):=\sum |\lambda_i|.$
Therefore with some constant $A^{\prime}$
%(we drop the $2\log n /n$ part
%in the definition of $\eps_n(f,\delta)$ for brevity),
\begin{eqnarray*}
&&
E_1\subseteq E_1^{\prime}:=\Bigl\{
\exists\delta\in (0,1) \ \exists f\in {\cal F}\
P_n\{f\leq \delta\}\leq
\Bigl(\frac{2\Delta(f)}{\delta}\Bigr)^{\frac{2\alpha}{2+\alpha}}
n^{-\frac{2}{2+\alpha}}\bigvee \frac{2\log n}{n}
\\
&&
{\rm and }\
P\{f\leq  \frac{\delta}{4}\}\geq A^{\prime}
\Bigl(
\Bigl(\frac{\Delta(f)}{\delta}\Bigr)^{\frac{2\alpha}{2+\alpha}}
n^{-\frac{2}{2+\alpha}}\bigvee \frac{2\log n}{n}
+\frac{t}{n}\Bigr)
\Bigr\}.
\end{eqnarray*}
Let us first consider the case when the class
$\cal H$ is uniformly bounded (say, by constant $1$).
One can observe that
${\cal F}'=\{f/\Delta(f) : f\in {\cal F}\}\subset
\{f\in {\rm conv}({\cal H}): \Delta(f)=1\}.$
For any function $f$ and any $\delta\geq \Delta(f),$
$P(f\leq\delta)=1,$ which means that on the event
$E_1^{\prime}$ one has to take into account
only values of $\delta\leq \Delta(f),$
or, equivalently, $\delta/\Delta(f)\leq 1.$
Therefore, a simple rescaling $\delta'=\delta/\Delta(f)<1$
shows that
$$
E_1^{\prime}=\Bigl\{
\exists\delta\in (0,1) \ \exists f\in {\cal F}'\
P_n\{f\leq \delta\}\leq
\Bigl(\frac{2}{\delta}\Bigr)^{\frac{2\alpha}{2+\alpha}}
n^{-\frac{2}{2+\alpha}}\bigvee \frac{2\log n}{n}\
{\rm and }\
$$
$$
P\{f\leq  \frac{\delta}{4}\}\geq A
\Bigl(
\Bigl(\frac{1}{\delta}\Bigr)^{\frac{2\alpha}{2+\alpha}}
n^{-\frac{2}{2+\alpha}}\bigvee \frac{2\log n}{n}
+\frac{t}{n}\Bigr)
\Bigr\}.
$$
As to the second condition on ${\cal F},$ in this case
$\Delta(f)=1$ for any $f$ by definition,
and the above equivalent representation of the event
$E_1^{\prime}$ holds automatically.

Let $\delta_j=2^{-j},\,j\geq 0.$
Theorem 2 (see also Example 1) and
a bound similar to (\ref{ecomp}) immediately imply that
for some $A$ and $B$
\begin{eqnarray*}
&&
{\mathbb P}\Bigl\{
\exists j\ \exists f\in {\cal F}'\
P_n\{f\leq \delta_j\}\leq
\Bigl(\frac{1}{\delta_j}\Bigr)^{\frac{2\alpha}{2+\alpha}}
n^{-\frac{2}{2+\alpha}}\bigvee \frac{2\log n}{n}\
{\rm and }\
\\
&&
P\{f\leq  \frac{\delta_j}{2}\}\geq A
\Bigl(
\Bigl(\frac{1}{\delta_j}\Bigr)^{\frac{2\alpha}{2+\alpha}}
n^{-\frac{2}{2+\alpha}}\bigvee \frac{2\log n}{n}
+\frac{t}{n}\Bigr)
\Bigr\}\leq
\\
&&
\leq\sum_{j\geq 0}B\exp\Bigl\{-\frac{1}{4}
\Bigl(\frac{\sqrt{n}}{\delta_j}\Bigr)
^{\frac{2\alpha}{2+\alpha}}\Bigr\}e^{-t/2}\leq
B'e^{-t/2}.
\end{eqnarray*}
The same argument as before yields
${\mathbb P}(E_1^{\prime})\leq Be^{-t/2}.$
Therefore, combining previuos bounds, we get
${\mathbb P}(E)\leq Be^{-t/4},$
which completes the proof of the theorem.

\qed

\section{Some experiments with learning algorithms}

In this section we present some results of the experiments
we conducted to test the ability of the new bounds
to predict the value of the generalization error of
combined classifiers. Unfortunately, the constants
in the bounds of Section 2 are not known.
More precisely, using the results of the recent work of
Massart \cite{Massart} one can calculate the constants involved
in the bounds, but their current values are rather large
and are way too far from being optimal. However, many important
learning algorithms (such as boosting and bagging)
that combine simple classifiers are
iterative in nature and it's important to see whether the
bounds allow one to predict the shape of the learning curves
(the dependence of the generalization error on the number
of iterations) correctly.
To this end, we just ignore the constants and use
in the experiments the quantities
$(n^{1-\gamma/2}\hat \delta_n(\gamma ;f)^{\gamma})^{-1}$
(see Example 1)
and $\eps_n(f;\hat\delta_n(f))$ (see Theorem
\ref{theorem:delta})\footnote{Actually, the quantity
$\eps_n(f;\hat\delta_n(f)/2)$ is involved in this bound;
but it's easy to see that it is within a constant
from $\eps_n(f;\hat\delta_n(f))$} instead of the upper bounds
we proved.
We will refer to these
quantities as the $\gamma$-bound and the $\Delta$-bound, respectively.
Incidentally, these quantities did provide
upper bounds on the generalization error (or on the test
error) in most of our experiments. This suggests that
the values of the constants involved in the bounds of Section 2
might actually be moderate (at least in the case when the bounds
are applied to several well known learning algorithms).

\subsection{Bagging and Boosting}

We begin by describing the experiments with two of
the most popular techniques
of combining the classifiers,
namely bagging \cite{Breiman:bagging}
and the Adaboost algorithm \cite{FS:Adaboost}.
In both of these methods,
there is an access to a learning algorithm
called a \emph{base learner}.
The base learner is given a training sample
$(X_i,Y_i),\ i=1,\dots ,n$ and it returns a classifier
$h$ from a base class $\HH$
that "approximately minimizes"
the empirical error $P_n\{yh(x)\leq 0\}$
(or properly weighted empirical error).

In the case of bagging, the base learner receives at each
iteration $t,$ $t=1,\dots ,T$ an independent bootstrap
sample $(\hat X_i^{(t)},\hat Y_i^{(t)}),\ i=1,\dots ,n$
and returns a classifier $h_t\in {\cal H}.$ The output
of bagging is the combined classifier $f:=T^{-1}\sum_{t=1}^T h_t$
(in other words, bagging makes a decision by majority vote).

%The combined classifier produced by
%Bagging (see algorithm \ref{alg:Bagging})
%has uniform weights.
%It is given by the majority vote of several
%base classifiers
%that are obtained by calling the base
%learner many times with different
%boostrap replicates of the original data set.
%That is, a new training set
%$\{X_{(i)},Y_{(i)}\}_{i=1}^n$
%is generated by drawing with
%replacement $n$ samples uniformly from the original data set.

%\begin{algorithm}
%\begin{algorithmic}
%\caption{Bagging}
%\label{alg:Bagging}
%\REQUIRE Sample $\{X_i,Y_i\}_{i=1}^n$, Base learning algorithm, Number of iterations $T$

%\FOR{$t=1$ to $T$}
%\STATE Generate boostrap sample $\{X_{(i)},Y_{(i)}\}_{i=1}^n$
%\STATE Find base classifier $h_t$ using Base learning algorithm and boostrap sample.
%\ENDFOR
%\STATE Output classifier $f(x)=\frac{1}{n}\sum_{t=1}^Th_t(x)$

%\end{algorithmic}
%\end{algorithm}

In the case of Adaboost, the algorithm assigns at the beginning
equal weights $D_1(i)=n^{-1},\ i=1,\dots ,n$
to all the training examples and
then updates the weights iteratively. Namely,
at $t$-th iteration ($t=1,\dots ,T$) the algorithm
calls the base learner that attempts to minimize
approximately the weighted training error
$$
\epsilon_{t}(h):=\sum_{i:h(X_i)\neq Y_i} D_t(i),\ h\in {\cal H}.
$$
The base learner returns a classifier $h_t\in {\cal H}$
and its weighted training error $\hat \epsilon_t:=\epsilon_t(h_t).$
The weights are then updated according to the formula
$$
D_{t+1}(i):=
\frac{D_t(i)}{Z_t}\bigl(1+(\beta_t-1)I_{\{h(X_i)=Y_i\}}\bigr),
$$
where $\beta_t:=\frac{\hat \epsilon_t}{1-\hat \epsilon_t}$
and $Z_t$ is the normalizing factor such that
$\sum_{i=1}^t D_{t+1}(i)=1.$ After $T$ iterations,
Adaboost outputs a combined classifier
$$
f:=\Bigl(\sum_{t=1}^T \log\frac{1}{\beta_t}\Bigr)^{-1}
\sum_{t=1}^T \log\frac{1}{\beta_t}h_t.
$$

%The base learners combined by Adaboost (see algorithm \ref{alg:Adaboost}) are produced by presenting the base learner  with weighted versions of the data set. The weights are initially uniform, and are adjusted at each step so that samples that are misclassified by the base classifier at the current step are given increased weight and samples that are classified correctly are given reduced weight. At each iteration, the base learner minimizes the empirical error with respect to the distribution given by
%he current values of the weights.

%begin{algorithm}
%begin{algorithmic}
%caption{AdaBoost}
%label{alg:Adaboost}
%REQUIRE Sample $\{X_i,Y_i\}_{i=1}^n$, Base learning algorithm, Number of iterations $T$
%STATE $D_1(i)=1/n$ for $i=1 \dots n $.
%FOR{$t=1$ to $T$}
%STATE Find base classifier $h_t$ using Base learning algorithm and distribution $D_t$.
%STATE $\epsilon_t=\sum_{i:h_t(X_i) \neq Y_i} D_t(i)$.
%STATE $\beta_t=\epsilon_t/(1-\epsilon_t)$
%FORALL{$i$}
%IF{$h_t(X_i)=Y_i$}
%STATE $D_{t+1}(i)=\frac{D_t(i)}{Z_t}\beta_t$
%ELSE
%STATE $D_{t+1}(i)=\frac{D_t(i)}{Z_t}$
%ENDIF
%STATE where $Z_t$ is a normalization factor so that $\sum_{i=1}^tD_{t+1}(i)=1.$
%ENDFOR
%ENDFOR
%STATE Output classifier $f(x)=\sum_{i=1}^T\log\frac{1}{\beta_t}h_t(x)$

%end{algorithmic}
%end{algorithm}

In all the experiments, we used
the set of indicator
functions\footnote{Actually, these functions are rescaled so that
they take values in $\{-1,1\}$}
of axis oriented hyperplanes
(also known as decision stumps) as base classifiers.
That is, $S:=\Reals^d$ and
$$
\HH=\left\{I_{\{\mathbf{x} \in \Reals^d: x_i \leq c\}},
c \in \Reals,\ i=1,\dots,d\right\}\cup
\left\{I_{\{\mathbf{x} \in \Reals^d: x_i \geq c\}},
c \in \Reals,\ i=1,\dots,d\right\},
$$
where $\mathbf{x}=(x_1,\dots,x_d)\in \Reals^d.$

\subsection{Experiments with real and simulated data}

We first describe the experiments with a "toy" problem
which is simple enough to allow one to compute exactly
the generalization error and other quantities such as
the $\gamma$-margins.
Namely, we consider a one dimensional classification problem
in which $S=[0,1]$ and, given a set (or a concept, using the
terminology of computer learning) $C_0\subset S$ which is
a finite union of disjoint intervals, the label $y$ is assigned to a point
$x\in S$ according to the rule $y=f_0(x),$ where $f_0$ is equal
to $+1$ on $C_0$ and to $-1$ on $S\setminus C_0.$
We refer to this problem as the {\it intervals problem.} 
Note that for the class of decision stumps we have
in this case $V(\HH)=2$
(since $\HH=\{I_{[0,b]}:b\in[0,1]\}\cup\{I_{[b,1]}:b\in[0,1]\}$),
and according to the results above the values of $\gamma $
in $[2/3,1)$ provide valid bounds on the generalization error
in terms of $\gamma$-margins.
In our experiments, the set $C_0$ was formed by
$20$ equally spaced intervals and we generated a uniformly
distributed on $[0,1]$ sample of size $1000.$
We ran Adaboost for 500
rounds (bagging does not work well for this problem),
and computed at each round the generalization error
of the combined classifier and
the quantity
$(n^{1-\gamma/2}\hat \delta_n(\gamma ;f)^{\gamma})^{-1}$
for different values of $\gamma$.

\begin{figure}
\begin{center}
\includegraphics[width=6cm]{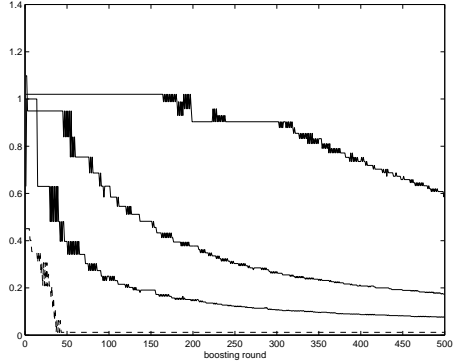}
\caption{Comparison of the generalization error (thicker line) with $(n^{1-\gamma/2}\hat \delta_n(\gamma ;f)^{\gamma})^{-1}$ for $\gamma=1,0.8$ and $2/3$ (thinner lines, top to bottom).}
\label{intervals_bounds1}
\end{center}
\end{figure}

\begin{figure}
\begin{center}
\includegraphics[width=6cm]{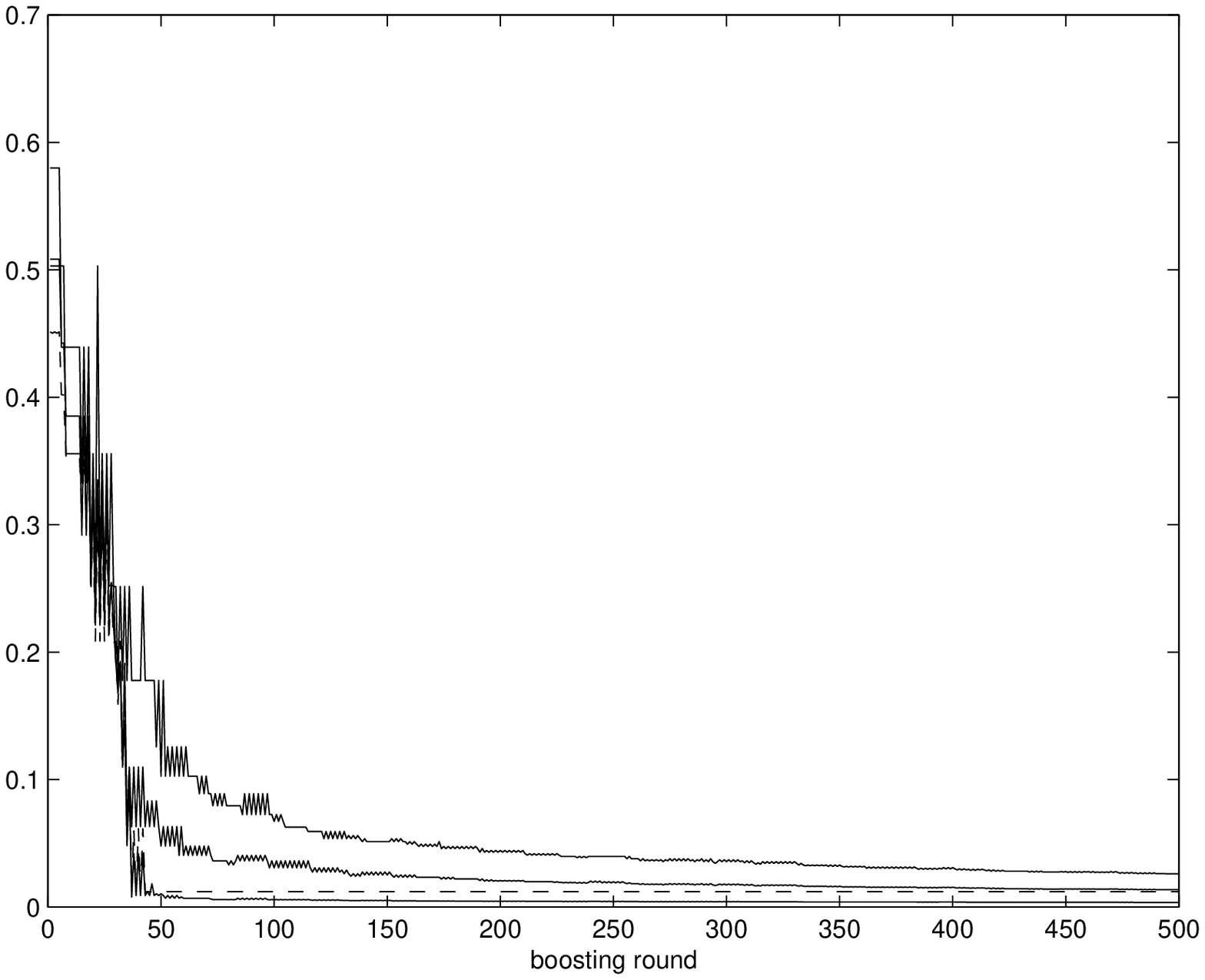}
\caption{Comparison of the generalization error (thicker line) with $(n^{1-\gamma/2}\hat \delta_n(\gamma ;f)^{\gamma})^{-1}$ for $\gamma=0.5,0.4$ and $0.2$ (thinner lines, top to bottom).}

\label{intervals_bounds2}
\end{center}
\end{figure}

\begin{figure}
\begin{center}
\includegraphics[width=6cm]{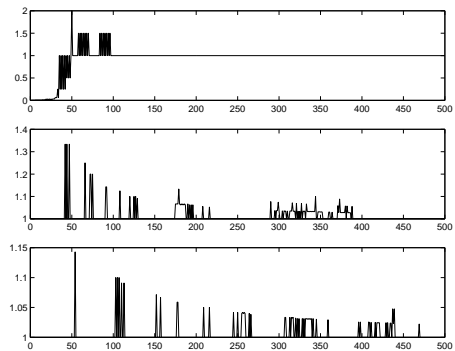}
\caption{Ratio $\hat{\delta_n}(\gamma ;f)/\delta_n(\gamma ;f)$ versus boosting round for
$\gamma=0.4,2/3,0.8$ (top to bottom)}
\label{ratio}
\end{center}
\end{figure}

In figure \ref{intervals_bounds1}
we plot the generalization error and the bounds for
$\gamma=1,0.8$ and $2/3$ against the iteration of Adaboost.
As expected, for $\gamma=1$
(which corresponds roughly to the bounds in \cite{SFBL})
the bound is very loose, and as $\gamma$ decreases,
the bound gets closer to the generalization error.
In figure \ref{intervals_bounds2} we show that by reducing
further the value of $\gamma$ we get a curve that is even closer
to the actual generalization error (although, for
$\gamma=0.2,$ it does not provide an upper bound
for some of the rounds of Adaboost). This
seems to support the conjecture that Adaboost actually generates
combined classifiers that belong to a subset of the convex hull
of $\HH$ with a smaller random entropy than of the whole
convex hull.
In figure \ref{ratio} we plot the ratio
$\hat{\delta_n}(\gamma ;f)/\delta_n(\gamma ;f)$
for $\gamma=0.4,2/3$ and $0.8$ against the boosting iteration.
We can see that the ratio is close to one in different examples
(for a small number of iterations of Adaboost in the first example,
the ratio is actually close to $0$)
indicating that the value of the constant $\bar A$ in
the bound (\ref{gamma1}) might be close to one
(at least, this seems to be true in the case of classifiers
produced by Adaboost for large sample sizes).

\begin{figure}
\begin{center}
\includegraphics[width=6cm]{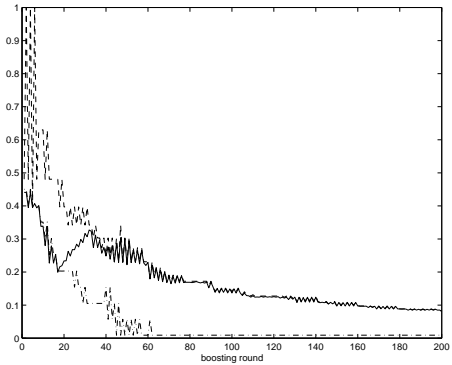}
\caption{Test error and bounds vs. number of classifiers
for the intervals problem for samples size of $1000$.
Test error (dot-dashed lines),
$\gamma$-margin bound with $\gamma=2/3$ (dashed lines),
and $\Delta$-bound (solid lines)}
\label{figure:bbounds_isel}
\end{center}
\end{figure}

In figure \ref{figure:bbounds_isel} we compare
the $\gamma$-bound and the $\Delta$-bound obtained for this problem for sample size of $1000$. We can see that the $\Delta$-bound has two regimes. In the first regime, the effect of the $\Delta$-dimension is dominant, and the bound tracks almost exactly the generalization error, giving a definite improvement over the $\gamma$-bound. In the second regime, the bound starts increasing until it reaches the curve of the $\gamma$-bound. This behavior can be explained
by examining the expression being minimized in the computation of the bound:

\begin{equation}\label{equation:two_terms}
\underbrace{
\frac{d(f;\Delta)}{n}\Bigl(\log\frac{1}{\delta}+
\log{\frac{ne^2}{d(f;\Delta)}}\Bigr)}_I+
\underbrace{\Bigl(\frac{\Delta}{\delta}\Bigr)^{\frac{2\alpha}{\alpha+2}}
n^{-\frac{2}{\alpha+2}}}_{II}
\end{equation}

It is easy to see that this expresion will be close to the
$\gamma$-bound when the second term is dominant,
and in fact, becomes the $\gamma$-bound when $\Delta=1$
(which, apparently, is the case in our experiments
when the number of classifiers in the convex combination becomes
large).

%This means that the effect of the $\Delta$-dimension
%on the bound is more important when the number of
%classifiers in the combination is relatively small,
%and as the number of classifiers increases,
%the margin term $II$ takes over.
%We observed that the rate of increase in the second
%regime is larger for a smaller sample size, which may be due
%to the factor $n$ in the denominator of $I$.

\begin{center}
\begin{figure}
\centering
\begin{tabular}{cc}
\hline
\textbf{Adaboost} &\textbf{Bagging}\\
\hline
\multicolumn{2}{l}{Twonorm}\\
\subfigure{\includegraphics[width=5cm,height=4cm]{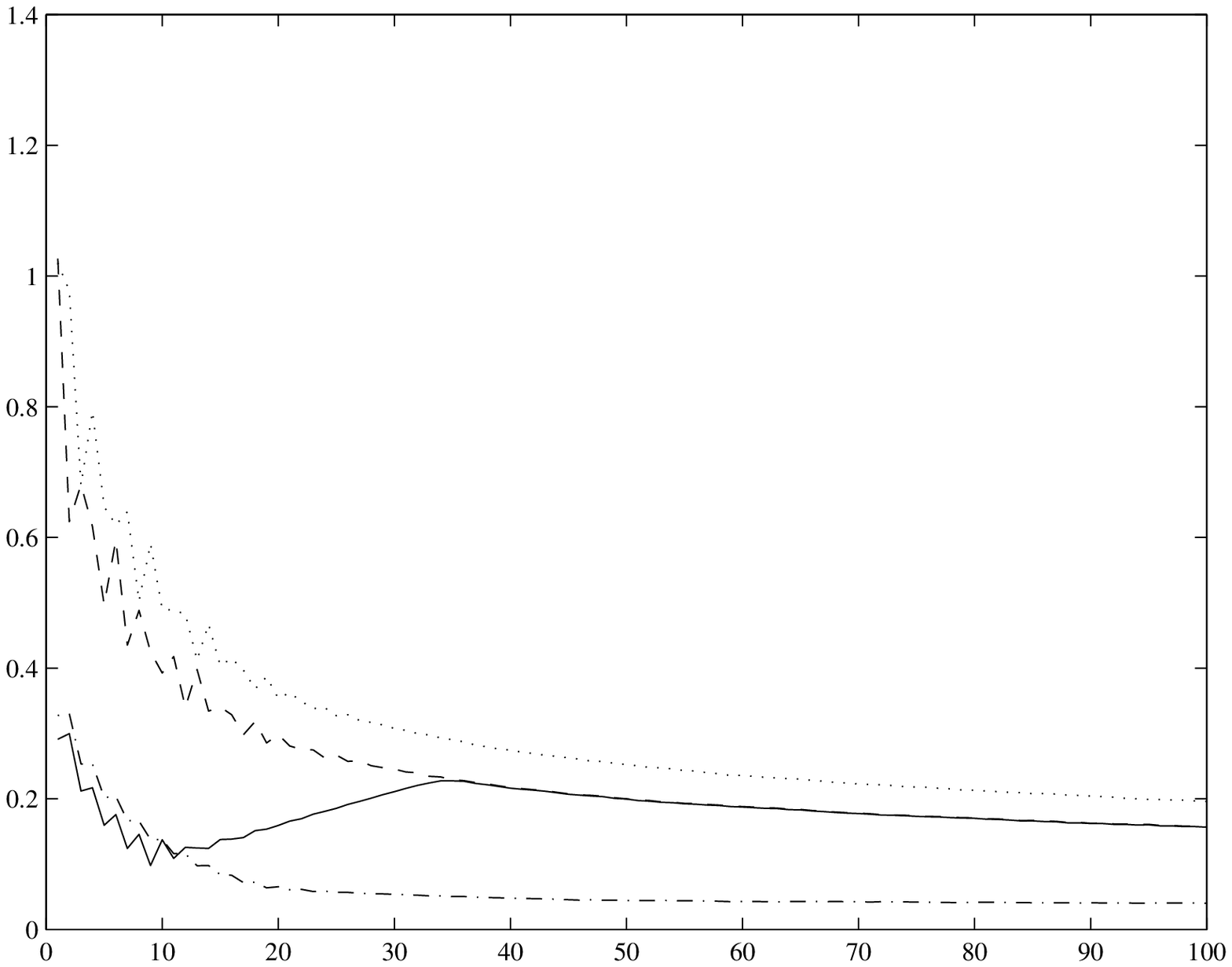}}&
\subfigure{\includegraphics[width=5cm,height=4cm]{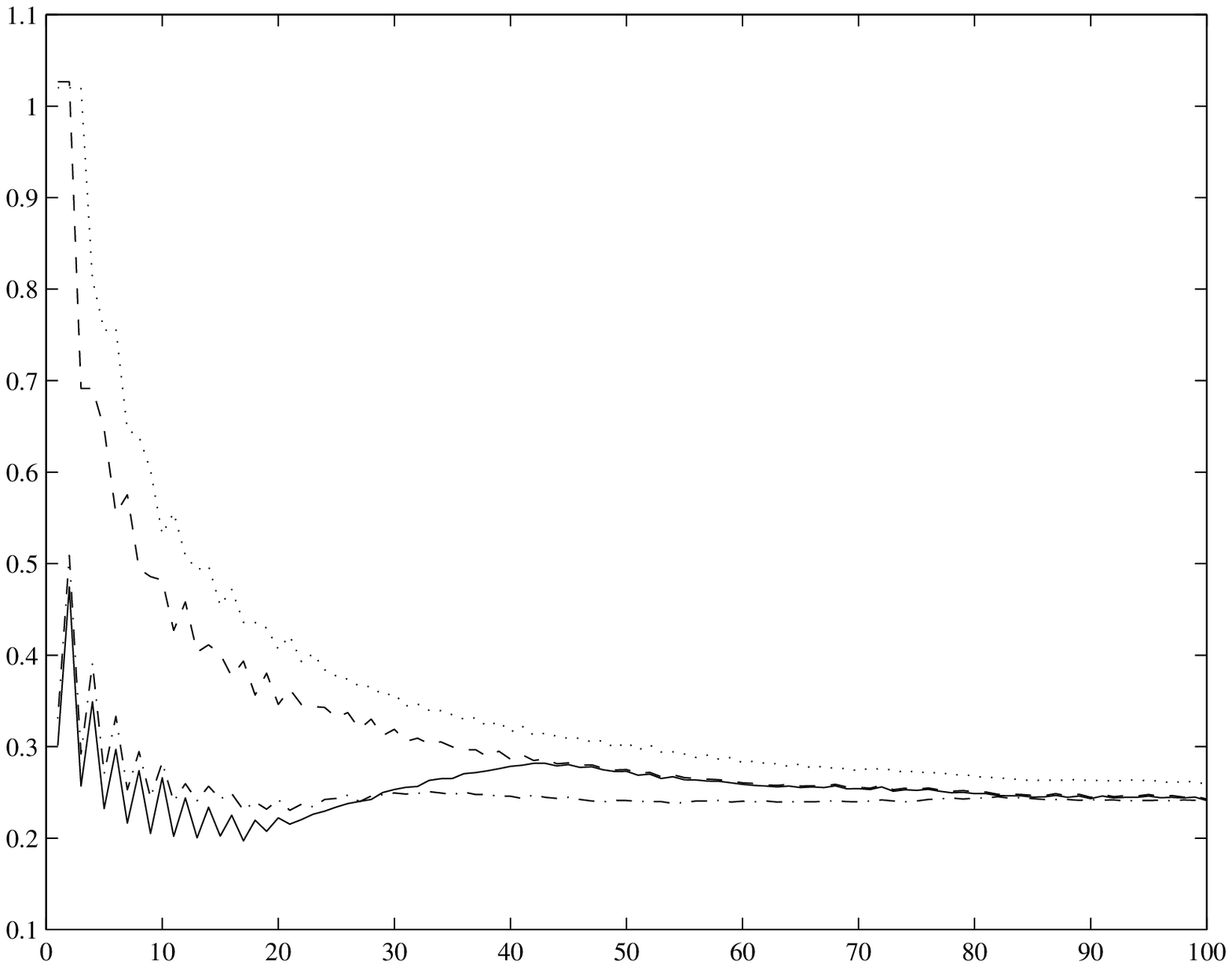}}\\
\hline
\multicolumn{2}{l}{King Rook vs. King Pawn}\\
\subfigure{\includegraphics[width=5cm,height=4cm]{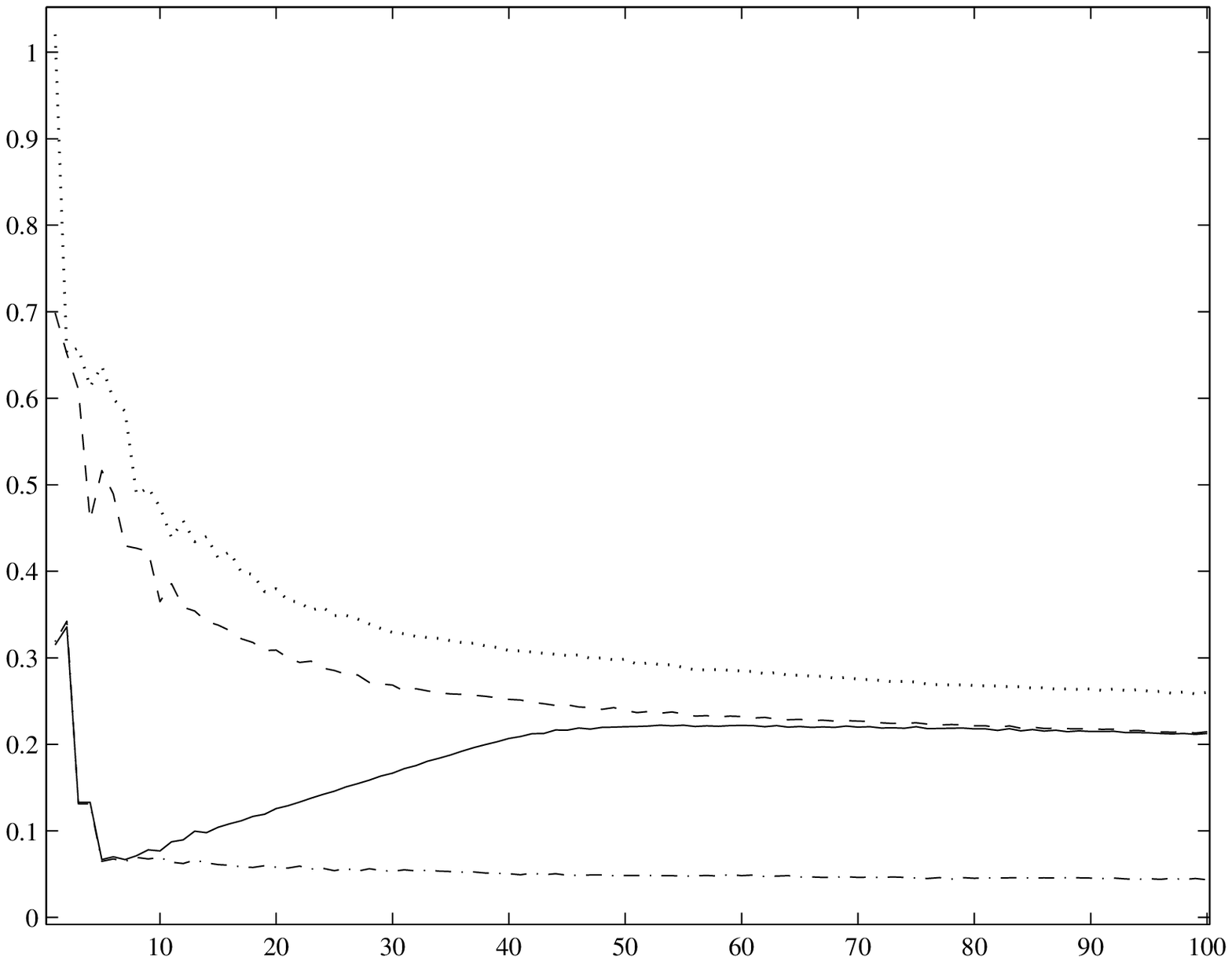}}&
\subfigure{\includegraphics[width=5cm,height=4cm]{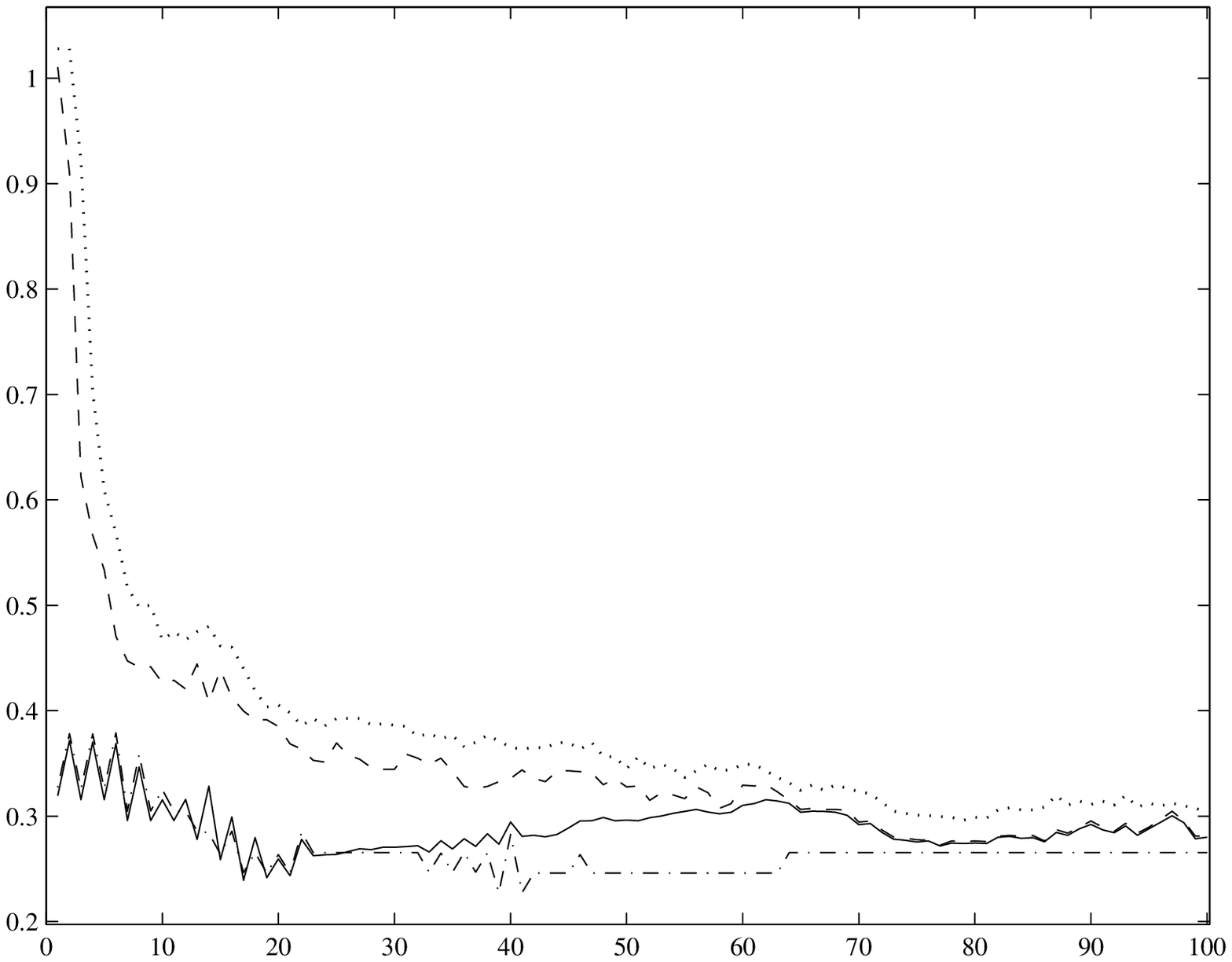}}\\
\hline
\end{tabular}
\caption{Test error and bounds vs. number of classifiers.
Test error (dot-dashed lines),
$\gamma$-margin bound with $\gamma=1$ (dotted lines), and
$\gamma=\gamma_{min}$ (dashed lines), and
the $\Delta$-bound (solid lines)}
\label{figure:bounds_real}
\end{figure}
\end{center}

We also computed the bounds for more
complex simulated data sets as well as for
real data sets in which the same type of behavior was observed.
We show the results for the so called Twonorm Data Set and the
King Rook vs. King Pawn Data Set (figure \ref{figure:bounds_real}),
which are well known examples in computer learning literature.
The Twonorm Data Set (taken from \cite{Breiman:arcing}) is a
simulated 20 dimensional data set in which positive and negative
training examples are drawn
from the multivariate normal distributions with unit covariance
matrix centered at $(2/\sqrt{20},\dots,2/\sqrt{20})$
and $(-2/\sqrt{20},\dots,-2/\sqrt{20}),$ respectively.
The King Rook vs. King Pawn Data Set is a real data set
from the UCI Irvine repository \cite{Blake+Merz:1998}).
It is a 36 dimensional data set with the sample size 3196.

As before, we used the decision stumps as base classifiers.
An upper bound on $V(\HH)$ for the class ${\cal H}$ of
decision stumps in $\Reals^d$ is given by the smallest
$n$ such that $2^{n-1}\geq (n-1)d+1$.
We computed the $\Delta$-bound and the $\gamma$-bounds
for $\gamma=1$ and for the smallest $\gamma $
allowed in Example 1 ($\gamma_{min}$).
For the Twonorm Data Set, we estimated the generalization error
by computing the empirical error on an
indepedently generated set of $20000$ observations.
For the King Rook vs. King Pawn Data Set,
we randomly selected $90\%$ of the data for training
and used the remaining $10\%$ to compute the test error.
The experiments were averaged over $10$ repetitions.

\subsection{Weighting and normalization}

It is apparent from the previous experiments that
the $\Delta$-bound explains well the behavior
of the generalization error for a small number of classifiers
in a convex combination, but for larger numbers of classifiers
it becomes close to $\gamma$-bound. Partially, it might be
related to the way the $\Delta$-dimension was defined.
In fact, the classifiers $h_t$ output
by the base learner at different iterations of Adaboost (or other
voting method of combining classifiers) can be close to each other
on the training examples (say, with respect
to the distance $d_{P_n,2}$).
Because of this, the $\Delta$-dimension
may very well overestimate the dimensionality of the combined
classifier
and more subtle definitions of dimension that take
into account such empirical closeness of different functions in
the convex combination are needed. The analysis of the proof of
Theorem 3 shows that the extension of our bounds to these more
subtle dimensions poses rather hard problems.

\begin{center}
\begin{figure}
\centering
\begin{tabular}{cc}
\hline
\textbf{Twonorm} &\textbf{King Rook vs. King Pawn}\\
\hline
\multicolumn{2}{c}{$\zeta=0.1$}\\
\subfigure{\includegraphics[width=4cm,height=3cm]{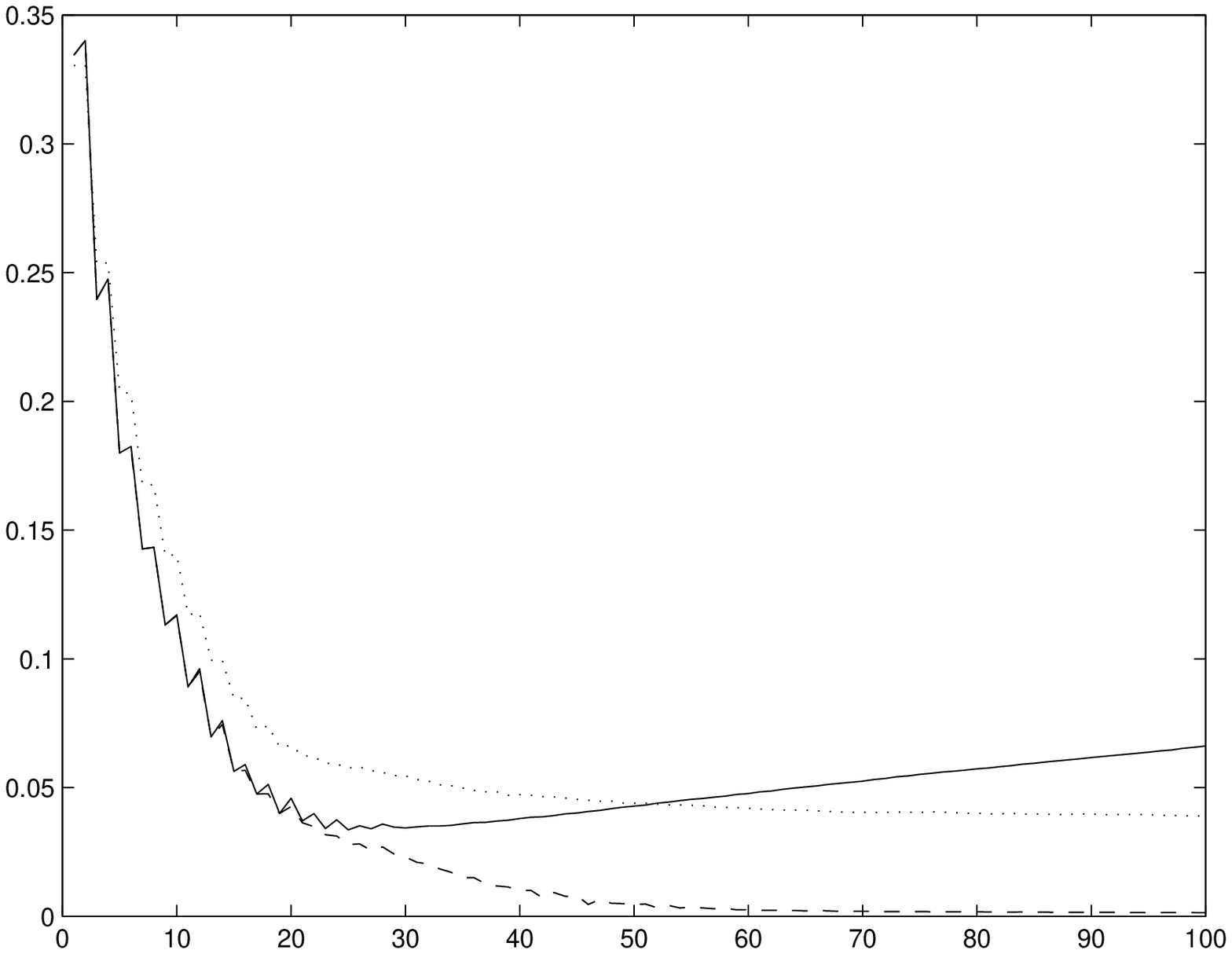}}&
\subfigure{\includegraphics[width=4cm,height=3cm]{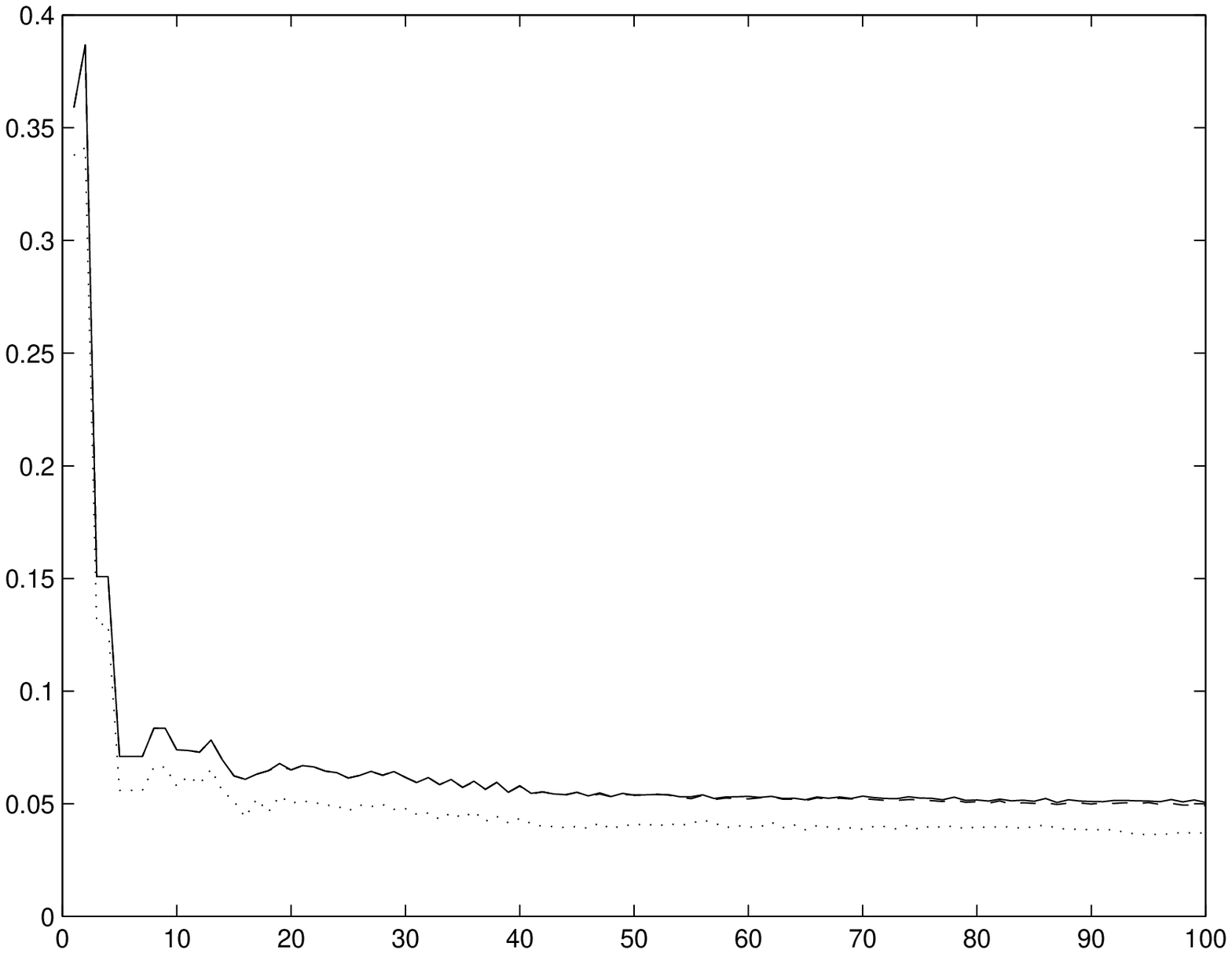}}\\

\hline
\multicolumn{2}{c}{$\zeta=0.4$}\\
\subfigure{\includegraphics[width=4cm,height=3cm]{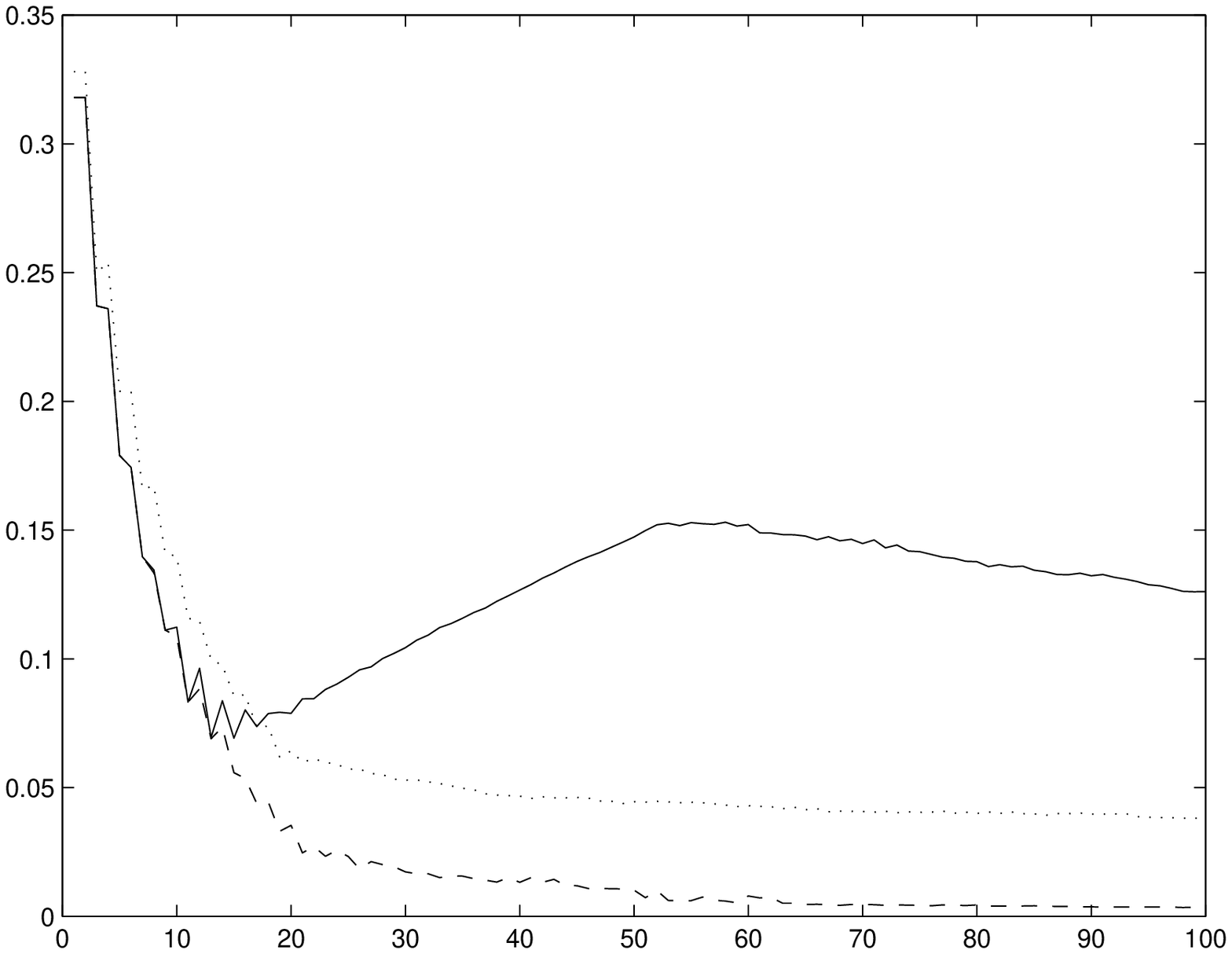}}&
\subfigure{\includegraphics[width=4cm,height=3cm]{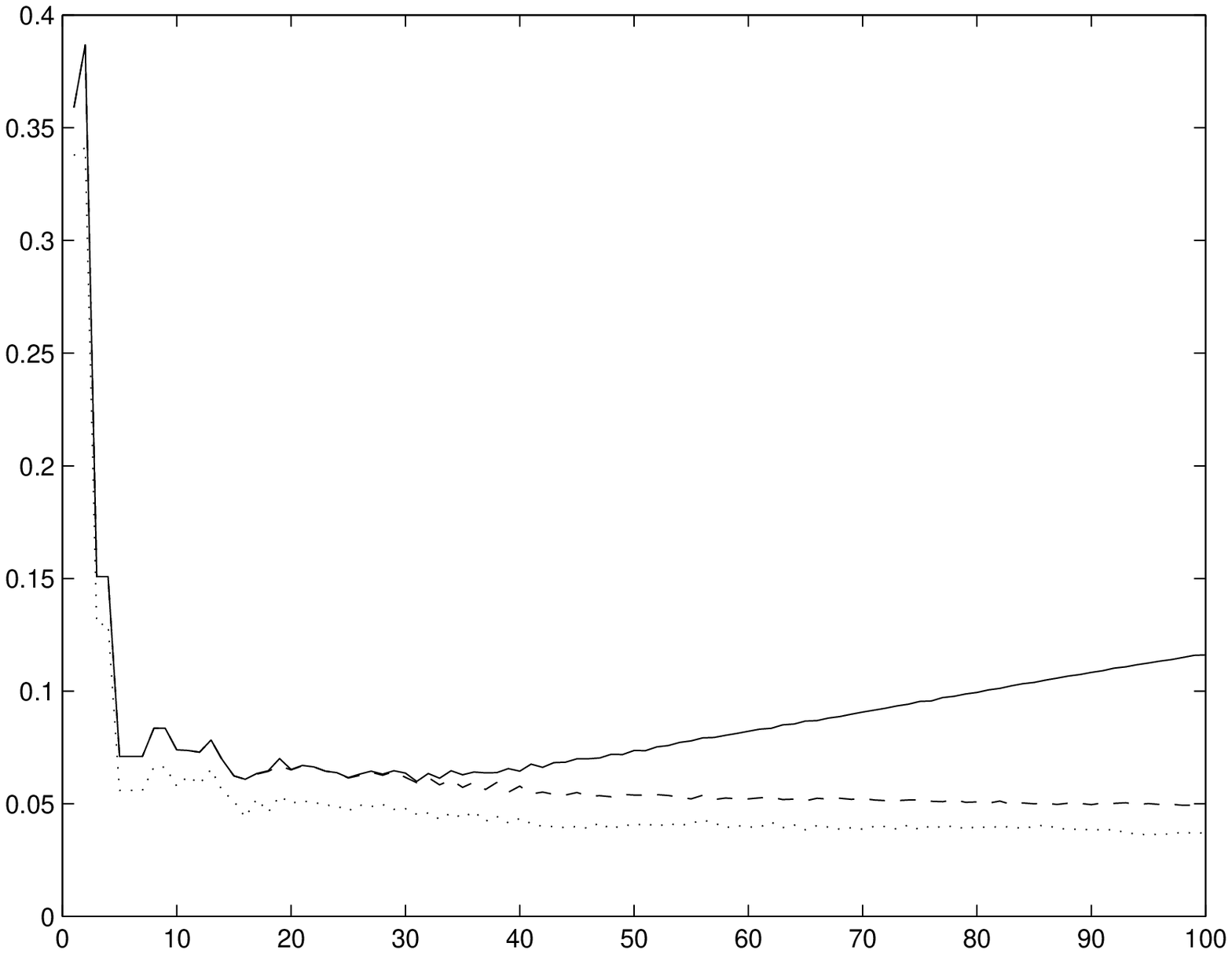}}\\
\hline
\multicolumn{2}{c}{$\zeta=0.9$}\\
\subfigure{\includegraphics[width=4cm,height=3cm]{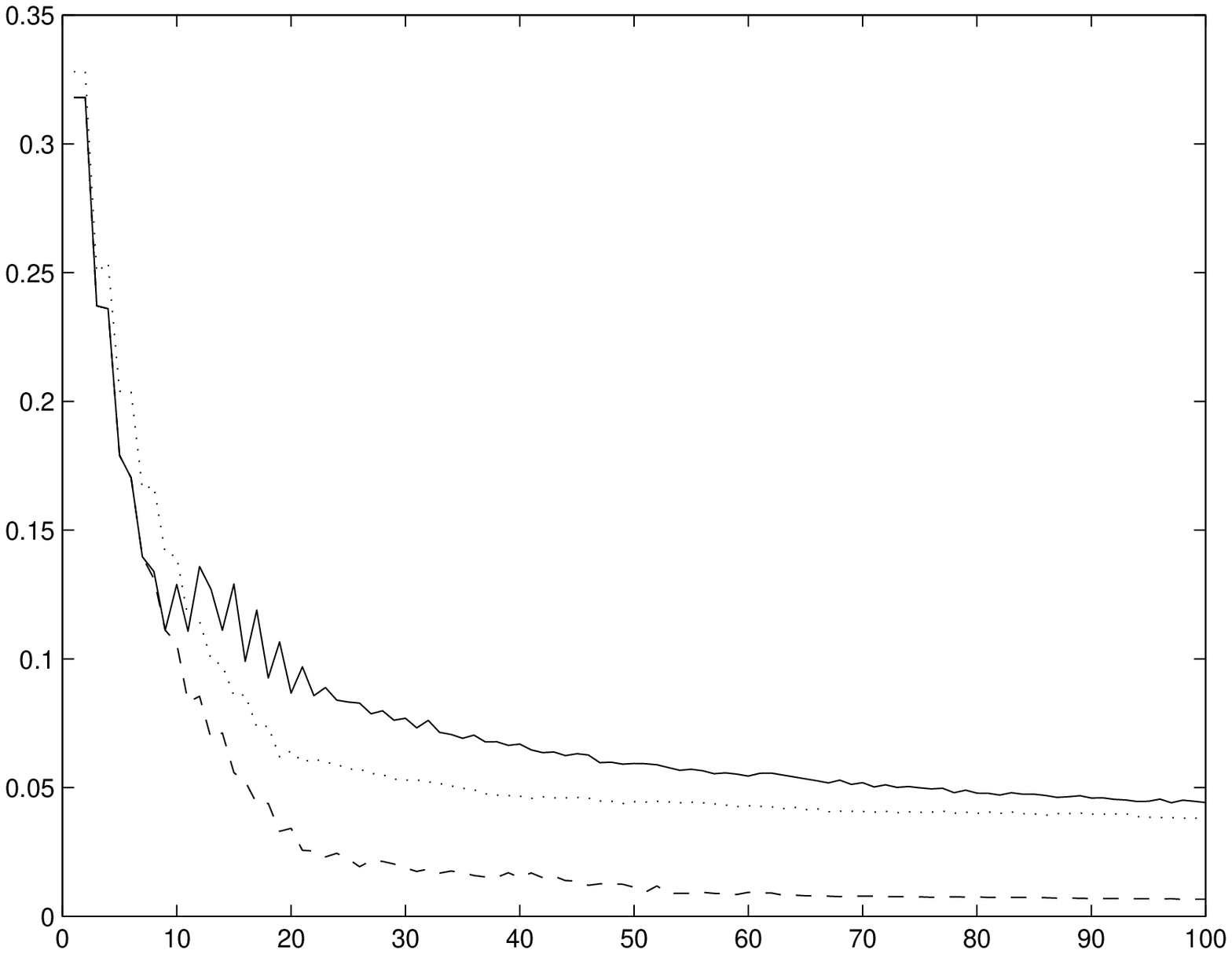}}&
\subfigure{\includegraphics[width=4cm,height=3cm]{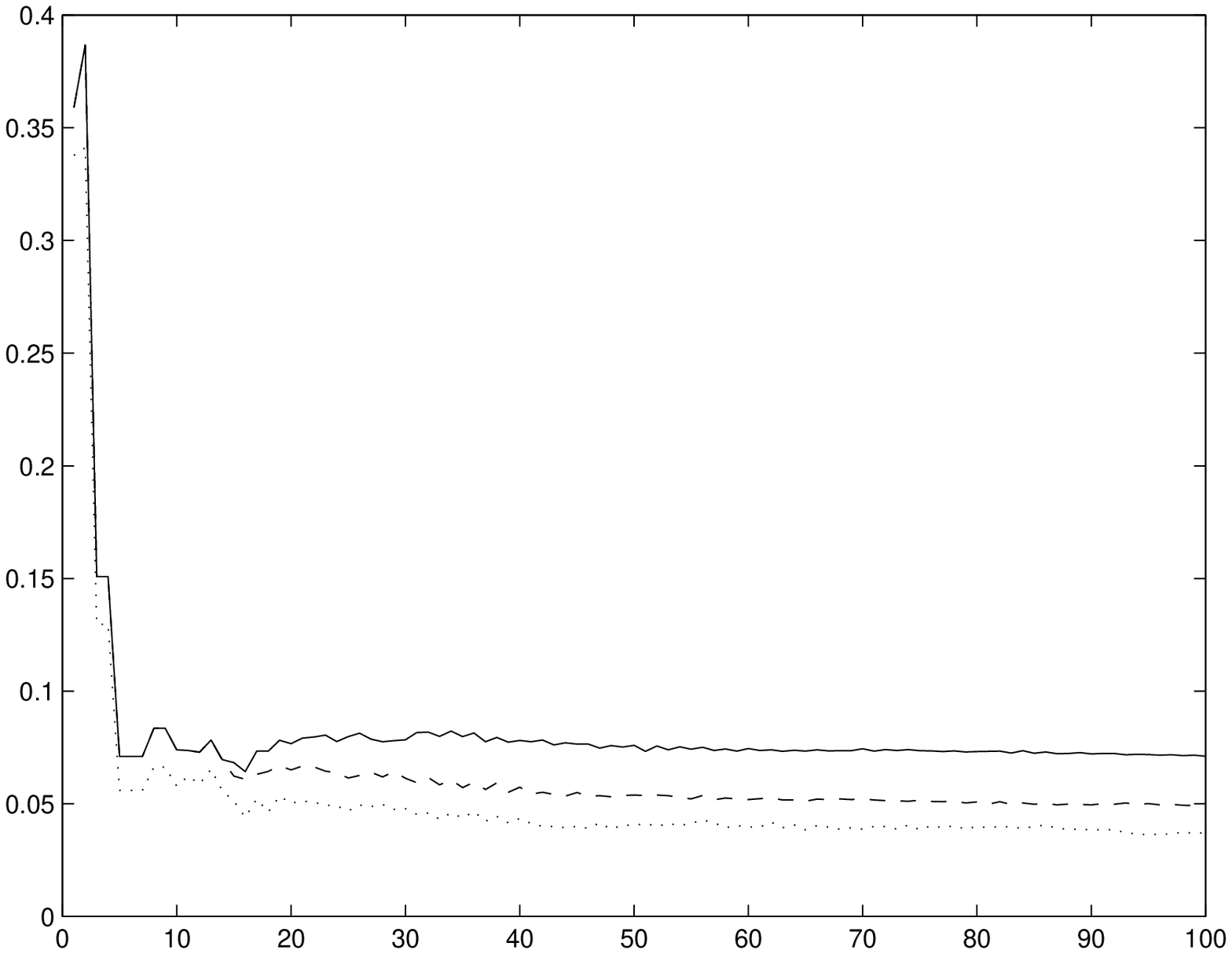}}\\
\hline
\end{tabular}
\caption{Bounds with weighting (solid line),
weighting and normalization (dashed line) and
test error (dotted line). In the bounds, $K=1.14$. }
\label{figure:w_and_n}
\end{figure}
\end{center}

%This is related to the increased value of the $I$ term
%in the expression (\ref{equation:two_terms})
%for a large number of classifiers,
%which causes the minimization procedure
%to select values of $\Delta$ for which the $\Delta$-bound
%becomes close to the $\gamma$-bound.
It might be also the case that the two terms in the expression
(\ref{equation:two_terms}) should be weighted in a certain
way in order to obtain a better bound.
The theoretical analysis of this problem is related
to determining sharp values of the constants involved in
the proof of Theorem 3 (which, in turn, is related to
the problem of optimizing the constants in Talagrand's concentration
and deviation inequalities for empirical processes that were used
in the proof). We performed some experiments in order to study
how such weighting influence the bound.
More precisely, given $\zeta \in [0,1]$ and $K>0,$
we defined
$$
\eps_{n,\zeta ,K}(f;\delta):=K\inf_{\Delta \in [0,1]}\Bigl[
\frac{\zeta d(f;\Delta)}{n}\Bigl(\log\frac{1}{\delta}+
\log{\frac{ne^2}{d(f;\Delta)}}\Bigr)+
(1-\zeta)\Bigl(\frac{\Delta}{\delta}\Bigr)^{\frac{2\alpha}{\alpha+2}}
n^{-\frac{2}{\alpha+2}}\Bigr]
$$
We also looked at a possibility of ``normalizing'' the value of
the $\Delta$-dimension in the bound
with respect to the total number of classifiers $T:$
$$
\tilde \eps_{n,\zeta,K}(f;\hat{\delta}_n(f)):=K\inf_{\Delta \in [0,1]}\Bigl[
\frac{\zeta d(f;\Delta)/T}{n}\Bigl(\log\frac{1}{\hat{\delta}_n(f)}+
\log{\frac{ne^2}{d(f;\Delta)}}\Bigr)+
(1-\zeta)\Bigl(\frac{\Delta}{\hat{\delta}_n(f)}\Bigr)^{\frac{2\alpha}{\alpha+2}}
n^{-\frac{2}{\alpha+2}}\Bigr].
$$
%The reason for this normalization is that when
%the number of classifiers increases,
%it may be possible to have a combined classifier
%with a large number of classifiers that have very small coefficients, making the $\Delta$-dimension large, even though most of the classifiers do not have any effect on the output of the combined classifier. This is the case for the Adaboost
%algorithm.
%We point out that we do not have at the moment any theoretical justification for this normalization.
%We conducted experiments to evaluate the behavior
%of the bound obtained when we
%incorporate weigthing and normalization in the computation.
We computed the bounds when weighting is used and when both
weighting and normalization are used.
%In the later case, we computed the following quantity:
We ran experiments for both simulated
and real data sets in which
we computed weighted and normalized bounds for values of
$\zeta=0.1,0.2,\dots 0.9$.
We show results for $\zeta=0.1,0.4$ and $0.9$ in figure
\ref{figure:w_and_n}.
%The value of the constant $K$ in these examples is  $1.14.$
%In most of the cases, we got an upper bound on the generalization
%error with this constant (although, in some examples, a larger
%value of $K$ was needed).

We found that weighting with a value of $\zeta=0.1$ gives
for most of the data sets a curve
that resembles rather closely the test error curve,
and does not present two different regimes as before.
%It should be mentioned that the values of $K$ and $\zeta$
%were determined for interval problem described at the beginning
%of Section 4.2. It happened that the same values also worked well
%in the experiments with other data sets.
When $\zeta$ increases (for example, when it becomes $0.4$)
the two-regime behavior becomes more noticeable,
although for $\zeta$ close to one the curves exhibit only
a small overshoot after which their shape is similar
to the shape of the test error curve.

When normalization is introduced, we get curves that are very
close to the test error curve for most of the data sets
(regardless of the value of parameter $\zeta$).
At the moment, we do not have any theoretical explanation
of these results.

\subsection{Towards algorithms balancing the dimensionality and the margins}

The connection between increasing the margins and reducing
the generalization error has led to the development
of several algorithms for designing and improving
combined classifiers based on optimizing
margin cost functions. The examples
include DOOM \cite{MBB}, DOOM2 \cite{MBBF},
DOOM-LP \cite{Lozano:Koltchinskii},
GeoLev \cite{Duffy:Helmbold},
and LP-Adaboost \cite{Grove:Schuurmans}.
The results in this paper motivate the development
of algorithms that take into account the
approximate dimensions of combined classifiers
along with their margins.

We discuss below the algorithm DOOM-LP, which was designed
to optimize a piecewise linear cost function
of the margins by solving a sequence of linear programs.
Incidentally, this algorithm also tends to reduce the
dimension of the combined classifier.
To describe the algorithm, define
$\varphi(u):=I_{(-\infty ,0]}(u)+(1-u)I_{(0,1]}(u)$
and let $\varphi_{\delta}(u):=\varphi(u/\delta).$
Let ${\cal H}$ be a base class and
${\cal F}:={\rm conv}({\cal H}).$
%Denote
%$$
%R_n({\cal H}):=\Bigl\|n^{-1}\sum_{i=1}^n \eps_i \delta_{X_i}\Bigr\|_{\cal H},
%$$
%where $\{\eps_i\}$ is a sequence of i.i.d. Rademacher random
%variables independent of the training data $(X_i,Y_i),\ i=1,\dots ,n.$
%The quantity $R_n({\cal H})$ is called the \it Rademacher complexity \rm
%of the class ${\cal H}.$
%If ${\cal H}$ is a VC-class,
%${\mathbb E}R_n({\cal H})\leq Cn^{-1/2}$ with a constant $C$
%depending on the VC-dimension of ${\cal H}.$
It was proved in Koltchinskii and Panchenko (2000)
that with probability at least $1-2\exp\{-2t^2\}$ the quantity
$$
\inf_{\delta\in [0,1]}\Bigl[P_n\varphi_{\delta}(yf(x))+
\frac{8}{\delta}{\mathbb E}\hat R_n({\cal H})+
\Bigl(\frac{\log\log_2(2\delta^{-1})}{n}\Bigr)^{1/2}\Bigr]
+\frac{t}{\sqrt{n}}
$$
is an upper bound on the generalization error $P\{yf(x)\leq 0\}$
of \it any \rm classifier $f\in {\cal F}.$
Recall that $\hat R_n({\cal H})$ is the Rademacher complexity
of the class ${\cal H}.$ If ${\cal H}$ is a VC-class, then
${\mathbb E}\hat R_n({\cal H})\leq Cn^{-1/2}$ with a constant $C$
depending on the VC-dimension of ${\cal H}.$
The idea of the algorithm DOOM-LP is to minimize the above
bound with respect to $f\in {\cal F}$ and $\delta\in [0,1]$
in order to find a classifier $\hat f$ with a reasonably
small generalization error. More precisely, the algorithms
receives a finite number of base classifiers $h_1,\dots, h_T$
along with their weights and attempts to redistribute the weights
in order to minimize the bound.
%DOOM-LP minimizes
%$P_n \varphi ({\frac{yf(x)}{\delta}})$
%for a large set of values of $\delta$ in $[0,1]$,
%computing the error of the resulting classifier
%on an independent validation set
%and selecting the set of weigths
%corresponding to the smallest error.

%Let $\varphi_{\delta}(m):=\varphi(m/\delta)$.

For a fixed value of $\delta $ and fixed classifiers
$h_1,\dots ,h_T,$ the minimization with respect
to $f=\sum_{k=1}^T w_k h_k\in {\cal F}$
consists of finding the weights $w_k,$
$\sum_{k=1}^T w_k=1,$
that minimize the following quantity:
\begin{equation}\label{eq:costf}
P_n\varphi_{\delta}(yf(x))=
\frac{1}{n}\sum_{i=1}^{n}\varphi_{\delta}
\left(Y_i\sum_{k=1}^{T}w_{k}h_k(X_i)\right).
\end{equation}
For a given combined classifier $f=\sum_{k=1}^T w_k h_k\in {\cal F},$
define sets $S_-,S_l,S_0$ as follows:
$$
S_-=\{i:Y_i f(X_i)\leq 0\},\
S_l=\{i: 0\leq Y_i f(X_i) \leq \delta\},\
S_0=\{i:Y_i f(X_i)\ge \delta\}.
$$
Finding the weight vector that "approximately minimizes"
$P_n\varphi_{\delta}(yf(x))$ for a fixed current partition
$(S_-,S_l,S_0)$ can be easily posed as a linear programming problem.
DOOM-LP searches for an approximate local
minimum of $P_n\varphi_{\delta}(yf(x))$
by solving this linear program and
moving to a neighboring partition by ``flipping'' the margins
that fall in the intersection of two of the sets $S_-,S_l,S_0$
from the set they currently belong to another one
in hope that with the constraints
determined by the new partition the objective function
can be reduced. The idea is similar in spirit to the
sweeping hinge algorithm proposed by Hush and Horn \cite{HH}.
The algorithm converges when the value of the minimum
in two neighboring partitions is the same (see algorithm
\ref{alg:DOOM-LP}).
We use the following notations in the description of the algorithm:
$b_k=-\sum_{i\in S_l}Y_ih_k(X_i)$ and $M_i=Y_i f(X_i),$
where  $f=\sum_{k} w_k h_k.$

\begin{algorithm}
\begin{algorithmic}
\caption{DOOM-LP}
\label{alg:DOOM-LP}
\REQUIRE Initial weight vector $\mathbf{w}$, margins $\{M_i\}_{i=1}^n$
\STATE \COMMENT {Initialize the partition}
\STATE $S_-=\{i:M_i\leq 0\}$
\STATE $S_l=\{i: 0\leq M_i \leq \delta\}$
\STATE $S_0=\{i: M_i \ge \delta\}$

\REPEAT
\STATE $C_{min}=\sum_{k=1}^Tb_kw_k$
\IF{$ |S_l| \geq 1 $}
\STATE \COMMENT {Compute optimal solution for a new partition}
\STATE $\mathbf{w}=\text{LPSolve}(\mathbf{w},S_-,S_l,S_0)$
\STATE Compute new margins $\{M_i\}_{i=1}^n$
\STATE \COMMENT {Update sets}
\STATE $S_-=S_- \cup \{i: i \in S_l, M_i=0\}-\{i: i \in S_-, M_{i}=0\}$
\STATE $S_l=S_l \cup \{i: i \in S_-, M_i=0\} \cup \{i: i \in S_0, M_{i}=\delta\}$
\STATE $\qquad -\{i: i \in S_l, M_i=0 \:\textrm{or}\: M_i=\delta\}$
\STATE $S_0=S_0 \cup \{i: i \in S_l, M_i=\delta\}-\{i: i \in S_0, M_{i}=\delta\}$
\STATE $C=\sum_{k=1}^Tb_kw_k$
\ELSE
\STATE Terminate and return current $\mathbf{w}$
\ENDIF
\UNTIL{$C\geq C_{min}$}
\end{algorithmic}
\end{algorithm}

If written in a standard form, the linear program solved
by DOOM-LP at each iteration involves $T+n+|S_l|+1$
variables ($T$ weights plus slack and surplus variables)
and $n+|S_l|+1$ equality constraints.
It follows from the basic results on linear programming
that if there is an optimal feasible solution and
the constraint matrix is full rank,
then there exists an optimal feasible solution with
at most $n+|S_l|+1$ non zero variables.
Furthermore, if the simplex method is used
to solve the linear program,
a solution of this type is allways found.
We have observed in experiments
that many of the variables that are set to zero in
the solution are weights and that DOOM-LP
tends to reduce the $\Delta$-dimension
of the classifier.

\begin{center}
\begin{figure}

\centering
\begin{tabular}{cc}

\subfigure[]{\includegraphics[width=5cm,height=4cm]{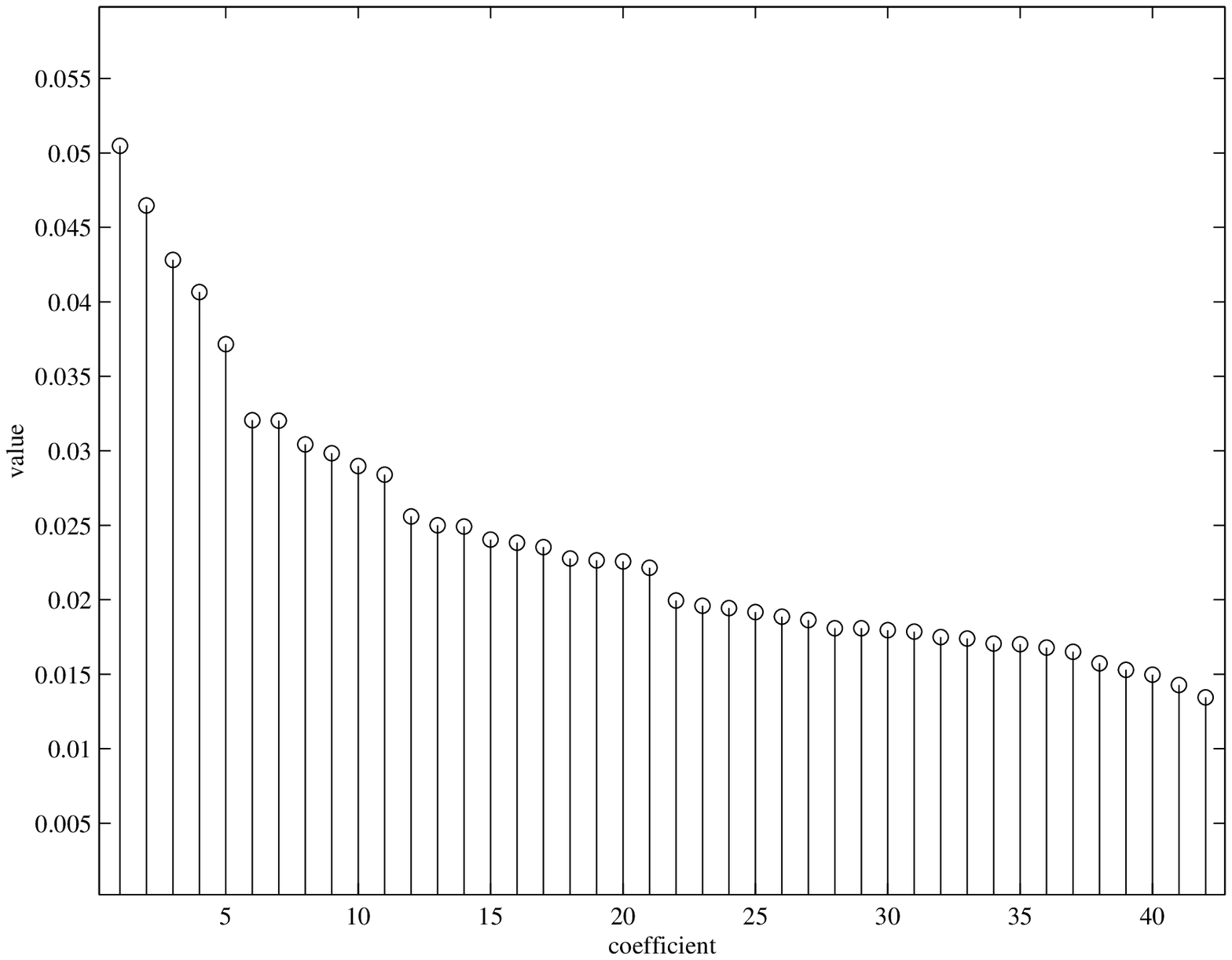}}&
\subfigure[]{\includegraphics[width=5cm,height=4cm]{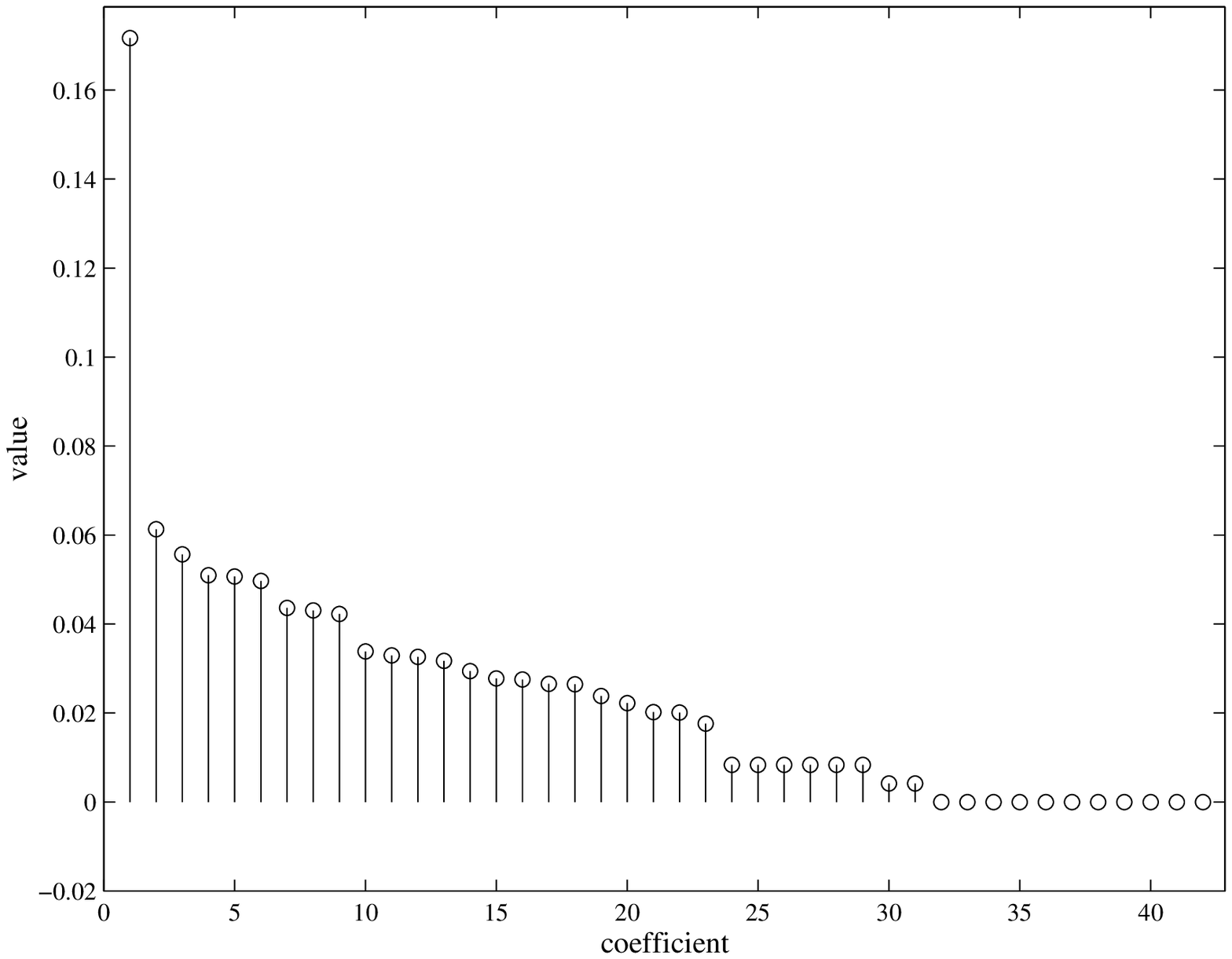}}\\
\subfigure[]{\includegraphics[width=5cm,height=4cm]{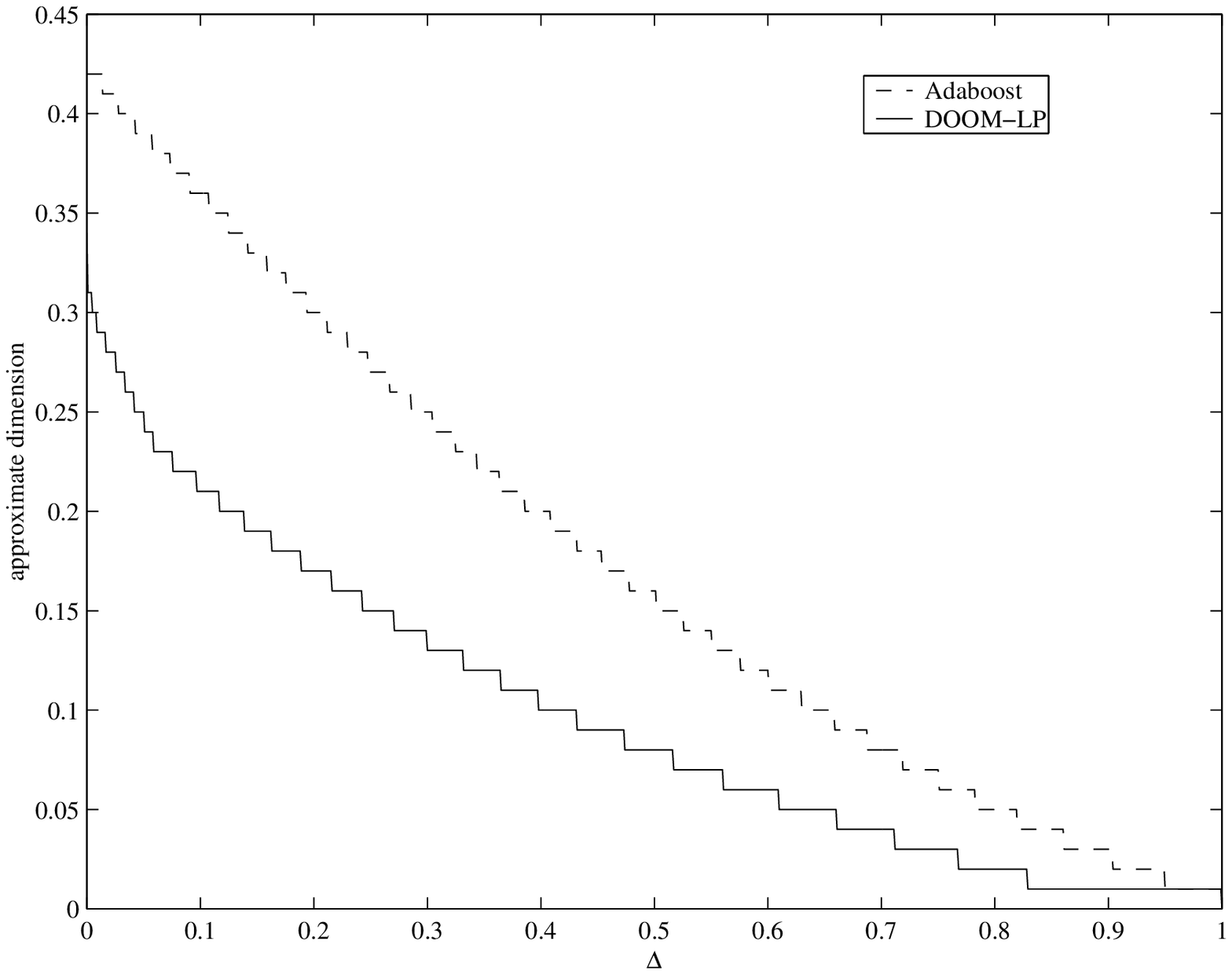}}\\
\end{tabular}
\caption{Results of running DOOM-LP on the classifier produced by Adaboost for the King Rook Vs. King Pawn data set. (a) Adaboost sorted coefficients, (b) DOOM-LP sorted coefficients, (c) Approximate $\Delta$-dimensions, (d) Cumulative margin distributions.}
\label{figure:kk}
\end{figure}
\end{center}

We have used DOOM-LP to improve the generalization error
of combined classifiers produced by Adaboost by redistributing
the weights of the base classifiers in a convex combination.
An example of dimensionality reduction by DOOM-LP
is illustrated in figure \ref{figure:kk}.

It might be interesting to design new algorithms
with explicit penalization for high
dimensionality in the optimization procedure.
For instance, assuming that the initial weights $w_t^{(0)}, t=1,\dots T$
are arranged in decreasing order,
one can add to the target function of linear program a
term $\sum_{t=1}^T a_t w_t,$ where $\{a_t, t\geq 1\}$ is
an increasing sequence of positive numbers. One can also consider entropy type penalties
of the form $\sum_{t=1}^T w_t \log\frac{1}{w_t}$ (in this case,
of course, the optimization is not a linear programming problem
any longer).

\hfill\break
Department of Mathematics and Statistics\hfill\break
The University of New Mexico\hfill\break
Albuquerque NM 87131-1141\hfill\break

\vskip 1mm

\hfill\break
Department of Electrical and Computer Engineering\hfill\break
The University of New Mexico\hfill\break
Albuquerque NM 87131\hfill\break


\begin{thebibliography}{99}
\bibitem{Anthony:Bartlett}
  M. Anthony, P. Bartlett (1999) 
Neural network learning : theoretical foundations,  
Cambridge University Press.
  


\bibitem{Baum:Haussler} E. Baum, D. Haussler (1989) 
What size net gives valid generalization?,  {\it Neural Computation}.
  




\bibitem{BBL2000} P. Bartlett, S. Boucheron, G. Lugosi (2000)  
Model Selection and Error Estimation. Preprint. 

\bibitem{Bartlett} P. Bartlett
 (1998)   
The sample complexity of pattern classification with neural networks: 
the size of the weights is more important than the size of the network, 
{\it IEEE Transactions on Information Theory}
  {\bf 44},
 525-536.
  

\bibitem{Bauer:Kohavi} E. Bauer, R. Kohavi
 (1999)
An Empirical Comparison of Voting Classification Algortithms: Bagging, Boosting and Variants,
{\it Machine Learning} {\bf 36},
 no. 1-2, 105-142.
  



\bibitem{Blum89} A. Blum, R. L. Rivest (1989)
Training a 3-node neural net is NP-Complete,  
Advances in Neural Information Processing Systems I,
 San Mateo, CA,
{494-501}.
  

\bibitem{Breiman:bagging} L. Breiman (1996)
Bagging Predictors,
  {\it Machine Learning},
 {\bf 26}, 123-140. 




\bibitem{Breiman:arcing} L. Breiman (1998) 
Arcing Classifiers, {\it The Annals of Statistics},
{\bf 26},
  no. 3,
  801-849.
 
 





\bibitem{Cherkauer}
  K. Cherkauer (1996)  
Human Expert-level performance on a scientific image analysis 
task by a system using combined artificial neural networks,   
Working notes of the AAAI Workshop on integrating multiple learned models,  
Eds: P. Chan, 15-21.  







\bibitem{Cortes:Vapnik} C. Cortes and V. Vapnik (1995)  
Support Vector Networks,
  {\it Machine Learning},
 {\bf 24}, 273-297.
 

\bibitem{Dietterich:Expboosting}
 T. Dietterich (2000) 
An Experimental Comparison of Three methods for constructing
 ensembles of decision trees: Bagging, Boosting and Randomization,  
{\it Machine Learning},
  {\bf 40}, 139-157.



\bibitem{Dietterich:survey} T. Dietterich (1997) 
Machine learning research: four current directions,  
{\it AI Magazine},
 {\bf 18,} no. 4, 97-136. 


\bibitem{Dietterich:Kong}
  T. Dietterich, E. B. Kong (1995)
Machine learning bias, statistical bias and statistical variance of decision tree algorithms, 
Department of Computer Science, Oregon State University.






\bibitem{Dietterich:Baikiri}
  T. Dietterich, G. Bakiri (1995)  
Solving multiclass learning problems via error-correcting output codes, 
{\it Journal of artificial intellingence research},
 {\bf 2}, 263-286. 







\bibitem{Drucker:Cortes}
  H. Drucker, C. Cortes (1996)  
Boosting decision trees,  
{Advances in Neural Information Processing systems 8},
 479-485.

\bibitem{Dudley}
  R.M. Dudley (1999)  
Uniform Central Limit Theorems,
 Cambridge University Press.








\bibitem{Duffy:Helmbold}
 N. Duffy, D. Helmbold (1999)
A geometric approach to leveraging weak learners,  
{Eurocolt99}.

\bibitem{Freund:Majority}
  Y. Freund (1995)  
Boosting a weak learning algorithm by majority, {\it Information and Computation}, 
{\bf 121}, no. 2, 256-285.

\bibitem{FS:expboosting}
  Y. Freund, R. Schapire (1996)  
Experiments with a new boosting algorithm,  
Machine Learning: Proceedings of the thirteenth international conference.

 



\bibitem{Grove:Schuurmans}
 A. Grove, D. Schuurmans (1998) 
Boosting in the limit: maximizing the margin of learned ensembles, 
Proceedings of the fifteenth national conference on Artificial intelligence.


\bibitem{Hush99}
 D. R. Hush (1999)
Training a Sigmoidal Node is Hard,
 {\it Neural Computation},
 {\bf 11},
1249-1260.






\bibitem{Kearns92}
  M. Kearns, R. E. Schapire, L. M. Sellie (1992)  
Towards efficient agnostic learning,
Proceedings of the 5th Annual Workshop on Computational Learning Theory,
New York, NY, ACM Press, 341-352.


  



\bibitem{Kolt_Control}
 V. Koltchinskii, C. Abdallah, M. Ariola, P. Dorato, D. Panchenko (1999)
Statistical learning control of uncertain systems: It is better than it seems,  
{Preprint, University of New Mexico}.
  

\bibitem{Kolt99}
 V. Koltchinskii (2001)    
Rademacher penalties and structural risk minimization,
{\it IEEE Transactions on Information Theory}, {\bf 47}, no. 5, 1902-1914.

\bibitem{KP00} V. Koltchinskii, D. Panchenko (2002)        
Rademacher processes and bounding the risk of function learning,    
High Dimensional Probability II.  



\bibitem{Kwok:Carter}
  S. Kwok,  C. Carter
 (1990)
Multiple decision trees,
 Uncertainty in Artificial Intelligence 4, 
Eds: R. Schachter and T. Levitt and L. Kannal and J. Lemmer,
Elsevier Science
, Amsterdam, 327-335.


\bibitem{Lozano:Koltchinskii}
  F. Lozano, V. Koltchinskii (2000)  
Direct optimization of simple cost functions of the margin,      
Preprint.
 

\bibitem{Massart}
 P. Massart (1998)   
About the constants in {T}alagrand's
  concentration inequalities for empirical processes.
 

\bibitem{Pitt_Valiant88}
  L. Pitt, L. Valiant (1988)  
Computational limitations of learning from examples, 
{\it Journal of the ACM},
 {\bf 35}, 965-984.


\bibitem{Pisier} G. Pisier (1981)   
Remarques sur un r\'esultat non publi\'e de {B.Maurey}.
S\'eminaire d'analyse Fonctionelle,
S\'eminaire d'analyse Fonctionelle, 1980-1981, Expos\'e No. 5, 
\'Ecole Polytechnique, Palaiseau.



\bibitem{Quinlan:boosting}
 J. R. Quinlan (1996) 
Bagging, boosting and C4.5,
 Fourteenth National Conference on Artificial Intelligence.
   








\bibitem{Schwenk:Bengio} H. Schwenk, Y. Bengio (1997)
Adaptive Boosting of Neural Networks for Character Recognition, 
D\'epartement d'informatique et recherche op\'erationnelle, Universit\'e de Montr\'eal. 
 

\bibitem{stone1}
 M. Stone (1974)   
Cross-validatory choice and assesment of statistical predictions, 
{\it Journal of the Royal Statistical Society},
 {\bf 36}, 111-147.


\bibitem{stone2}
  M. Stone (1977)    
Asymptotics for and against cross-validation,  
{\it Biometrika},
 {\bf 64}, no. 1, 29-36.


\bibitem{Talagrand1}  M. Talagrand (1996)  
A new look at independence,
 {\it Annals of Probability}, 
{\bf 24},  1-34.


\bibitem{Talagrand2}
  M. Talagrand (1996)   
New concentration inequalities in product spaces,  
{\it Invent. Math.},
  {\bf 126}, 505-563.

\bibitem{VandW}  A. W. van der Vaart, J. Wellner (1996) 
Weak Convergence of Empirical Processes With Applications to Statistics,
 Springer.
   


\bibitem{VAP98}  V. Vapnik (1998)  
Statistical Learning Theory,
 John Wiley and Sons, Inc.






\bibitem{Vidyasagar}
  M. Vidyasagar (1997) A theory of learning and generalization,
 Springer-Verlag.


\bibitem{DVGL}
  L. Devroye, L. Gy\"orfi, G. Lugosi (1996)
A Probabilistic Theory of Pattern Recognition,
  Springer-Verlag.


\bibitem{FS:Adaboost}
  Y. Freund, R. Schapire (1997)  
 A decision-theoretic generalization of on-line learning and an application to boosting,
{\it Journal of Computer and System Sciences},
 {\bf 55}, no. 1, 119-139. 



\bibitem{HH}
  D. Hush and B. Horne (1998)  
Efficient algorithms for function approximation with piecewise linear sigmoids,
{\it IEEE Transactions on Neural Networks},
 {\bf 9}, no. 6, 1129-1141. 


\bibitem{KP}
 V. Koltchinskii, D. Panchenko (2002)  
Empirical margin distribution and bounding the generalization error of
combined classifiers,
 {\it Annals of Statistics},
 {\bf 30}, no. 1.

\bibitem{Kol2001}  V. Koltchinskii (2001) 
Bounds on margin distributions in learning problems,  
Preprint.

\bibitem{KPL}
  V. Koltchinskii, D. Panchenko, F. Lozano (2000) 
Bounding the generalization error of neural networks and combined classifiers,  
{Advances in Neural Information Processing Systems 13}.

 





\bibitem{KPL2} V. Koltchinskii,  D. Panchenko, F. Lozano
 (2001) 
Bounding the generalization error of convex combinations of classifiers: 
balancing the
 dimensionality and the margins. Preprint. 
 

\bibitem{SFBL} R.E. Schapire, Y. Freund, P. Bartlett, W.S. Lee (1998)  
Boosting the margin : A new explanation for the efectiveness of voting methods,
{\it Annals of Statistics},
 {\bf 26}, no. 5, 1651-1687.

\bibitem{MBB}
  L. Mason, P. Bartlett, J. Baxter (2000) 
Improved generalization through explicit optimization of margins,
{\it Machine Learning}, {\bf 38},
  no. 3,
 243-255.


\bibitem{MBBF}
 L. Mason, J. Baxter, P. Bartlett, M. Frean (1999)
Advances in large margin classifiers,
 MIT Press.
 
\bibitem{Blake+Merz:1998}  C.L. Blake, C.J. Merz (1998)  
{UCI} Repository of machine learning databases,   
URL: http://www.ics.uci.edu/$\sim$mlearn/MLRepository.html,
   
University of California, Irvine, Dept. of Information and Computer Sciences.

\bibitem{Mitchell}
 M. Mitchell (1996)  
An Introduction to Genetic Algorithms,
  MIT Press.

 





\bibitem{DeJong}
  K.A. DeJong, M. Spears, D. Gordon (1993)  
Using Genetic Algorithms for Concept Learning,  
{\it Machine Learning},
  {\bf 13}, 161-188.






\bibitem{Freund99}  Y. Freund (1998) Self Bounding Learning Algorithms,  
Proceedings of the Eleventh Annual Conference on Computational Learning Theory.


  

\bibitem{Holland75}
  J.H. Holland (1975)  
Adaptation in Natural and Artificial Systems, 
University of Michigan Press.




\bibitem{Maniezzo}  V. Maniezzo (1994) 
Genetic Evolution of the Topology and Weight Distribution of Neural Networks,
  
{\it IEEE Transactions on Neural Networks},
 {\bf 5}, no. 1, 39-53. 

 
\bibitem{Montana}  D.J. Montana, Davis L. (1989) 
Training Feeforward Neural Networks using Genetic Algorithms, 
Proceedings of the International Joint Conference on Artificial Intellingence,
Morgan Kaufmann.





\bibitem{Valiant:pac}
  L. Valiant (1984)   
A theory of the learnable,
 {\it Communications of the ACM},  
{\bf 27}, no. 11, 1134-1142.


\bibitem{Yao}
  X. Yao (1993) 
A review of evolutionary artificial neural networks,  
{\it International Journal of Intelligent Systems},
 {\bf 8}, no. 4, 539-567. 




\end{thebibliography}
\end{document}